\documentclass[a4paper, freqn]{amsart}

\usepackage[T2A,T1]{fontenc}  % T2A for Cyrillic support
\usepackage[russian, french, british]{babel}
\usepackage{mathtools}
\usepackage{bm}
\usepackage[hidelinks]{hyperref}
\usepackage{amssymb, amsmath}
\usepackage{comment}
\usepackage{url}

\renewcommand{\phi}{\varphi}

\DeclarePairedDelimiterX{\Set}[2]{\{}{\}}{ #1 \mathchoice{\:}{\:}{\,}{\,}\delimsize\vert\allowbreak\mathchoice{\:}{\:}{\,}{\,}\mathopen{} #2 }
\DeclarePairedDelimiterX{\Seq}[2]{\langle}{\rangle}{ #1 \mathchoice{\:}{\:}{\,}{\,}\delimsize\vert\allowbreak\mathchoice{\:}{\:}{\,}{\,}\mathopen{} #2 }

\newcommand{\CC}{\mathcal C}
\newcommand{\LL}{\mathcal L}
\newcommand{\PP}{\mathcal P}

\newcommand{\Aa}{\mathbf A}
\newcommand{\Bb}{\mathbf B}

\newcommand{\cfr}{\mathfrak c}

\theoremstyle{plain}
\newtheorem{theorem}{Theorem}[section]
\newtheorem{observation}[theorem]{Observation}
\newtheorem{proposition}[theorem]{Proposition}
\newtheorem{lemma}[theorem]{Lemma}
\newtheorem{corollary}[theorem]{Corollary}

\theoremstyle{definition}
\newtheorem{definition}[theorem]{Definition}

\theoremstyle{remark}
\newtheorem{remark}[theorem]{Remark}
\newtheorem{notation}[theorem]{Notation}
\newtheorem{example}[theorem]{Example}

\newcommand{\satisfies}{\vDash}
\newcommand{\rest}{\upharpoonright}
\newcommand{\wl}{{\mathsf {wl}}}

\begin{document}
\title[MSO Theory of the Reals Order with the Borel Sets Quantifiers]{Monadic Second Order Logic of the Real Order with the Set Quantifiers ranging over the Borel Sets}

\author{Mirna D{\v z}amonja}
\address{base-name Logique Consult, 7 Impasse Charti\`{e}re, 75005 Paris, France}
\email{logique.consult@gmail.com}
\urladdr{https://www.logiqueconsult.eu}

\thanks{I thank Joan Bagaria for the invitation to visit him at the University of Barcelona in April 2025. The conversation we had gave me a crucial insight. Namely, I showed Joan the idea of well-levelled orders and he proposed that this could probably be used in descriptive set theory. I went to work on this and connected it with Rabin's work and solving the problem of MSO on the Borel sets, about which I had been thinking about on and off for a while. The idea of writing a joint paper on this did not converge, due to Joan's other projects. I thank Andreas Blass, Ohad Drucker, Yair Hayut, Jouko Väänänen and Andrés Villaveces on very useful comments on the preprint of the paper. In addition, I thank both the Set Theory Seminar and the Logic Seminar  at the Hebrew University of Jerusalem for the active comments of the participants and for giving me an opportunity to present this work at its and mine alma mater.}

\thanks{{\em Declaration on the use of artificial intelligence}: I use Grok (by xAI) as a research assistant that does everything that a traditional secretary could do, although I do my own typing. In particular, Grok can search for a reference, help with a {\LaTeX} macro or with English. I have not knowingly used Grok as a mathematics assistant, although it is likely capable of proof checking. All mathematical content, proofs and text were authored and written by me and I take the full responsibility for the paper's accuracy and integrity.}

\maketitle

\begin{abstract} A celebrated 1969 theorem of Michael Rabin is that the monadic second order (MSO) theory of the real order where the monadic quantifier is allowed only to range over the closed or $F_\sigma$-sets, is decidable. In 1975 Saharon Shelah \cite[Conjecture 7A]{Shelahmonadic}
conjectured that if the monadic quantifier is allowed to range over the Borel subsets of the reals, the resulting MSO theory is still decidable. We confirm this conjecture. In fact, the conjecture can be understood in the weak form, in the language where there is a special symbol for each level of the Borel hierarchy, or in a strong form where the language simply has a symbol for Borel sets. We confirm both versions of the conjecture, the former one by interpreting it in S2S and the latter one by in addition using Büchi's theorem that MSO$(\omega_1,<)$ is decidable.
\end{abstract}

\section{Introduction} Monadic Second Order logic (MSO) is defined like the first order logic (FO), but with an addition of set variables and the additional atomic formulas of the form $x\in {\mathbf A}$ and $(\exists  {\mathbf A})\psi({\mathbf A})$ for a set-quantifier-free $\psi({\mathbf A})$. This logic is a small fragment of the Second Order Logic (SO), yet it has proven sufficiently rich to express various useful mathematical statements. The results of this paper show, for example, that for any given $\alpha<\omega_1$ the Determinacy of the $\alpha$-level of the Borel hierarchy is expressible as a sentence in the MSO theory S2S of the
Cantor tree with the two successor functions (Corollary \ref{th:effectiveBorel}). On the other hand, MSO is close enough to FO to have many properties that let us work with the logic in practice, in particular many instances of decidability. The decidability of various logics is a question stemming directly from Hilbert's programme, most famously through Hilbert's 10th Problem and Yuri Matijashevich's 1970 solution to it \cite{zbMATH03336816}. The full SO is most often undecidable. In the case of SO arithmetic, so the SO theory of $({\mathbb N}, +,\cdot, <, 0,1)$, the undecidability was proven by Alonzo Church in \cite{88d09d93-3574-3ee8-8658-ea732b3f5a1b}. In contrast, the MSO theory of various structures is decidable, including the MSO theory of $({\mathbb N}, +,\cdot, <, 0,1)$, which was famously proved decidable by J. Richard Büchi in \cite{Buchi}. 
Büchi uses an important method, which is to associate to an MSO theory the runs of an automaton. This method was introduced by Boris (Boaz) Trakhtenbrot in \cite{trakhtenbrot1958english}, and in fact, a posteriori, one could say that the most important property of MSO is its closeness to the theory of automata. This makes MSO a logic much studied in theoretical computer science. The survey article \cite{GurevichMSO} by Yuri Gurevich presents many properties of the MSO logic. 

A celebrated theorem of Michael Rabin in \cite[Th 2.7]{Rabin} vastly extended Büchi's theorem in proving that the MSO theory of the real order where the monadic quantifier is allowed only to range over closed or $F_\sigma$-sets, is decidable. The limits of generalisations were shown by Saharon Shelah
in an important 1975 paper \cite{Shelahmonadic}. Shelah proved that under the assumption of the continuum hypothesis (CH), if the monadic quantifier is allowed to range over all subsets of the reals, the resulting MSO theory is undecidable (the assumption of CH was removed by Gurevich and Shelah in \cite{GUREVICH1982179}). Shelah conjectured \cite[Conjecture 7A]{Shelahmonadic} that if we allow the monadic quantifier to range over the Borel subsets of the reals, the resulting theory is decidable.

In Theorem \ref{th:main} and Corollary \ref{th:strongmain} we confirm Shelah's conjecture. The MSO theory of the real order where the monadic quantifier is allowed to range over the Borel subsets of the reals, is decidable.

Our methods are similar and inspired by those used by Rabin in the proof of \cite[Th 2.7]{Rabin}  to code closed and $F_\sigma$-sets in S2S, except that we have developed an additional coding technique. It has allowed us to code in S2S Borel sets of arbitrary complexity.  
That technique is a specific instance of a more general idea of {\em well-levelled orders} that this paper introduces.

The paper is organised as follows. In \S\ref{sec:def} we define well-levelled orders. They are a class of well-founded partial orders and every well-order and every tree is well-levelled. We present the basic properties of well-levelled partial orders  and show some examples, to the extent that will be useful for the main result and slightly more, as a reference. In \S\ref{sec:realcoding} we show how to use 
these ideas to code Borel subsets of the reals. 
In \S\ref{sec:automates} we show how to use this coding, combined with the ideas from the original proof by Rabin, for the main decidability result. The main results are presented in \S\ref{sec:Borel} in Theorem \ref{th:main} and Corollary \ref{th:strongmain}. For the stronger version of the main result, which is 
Corollary \ref{th:strongmain}, we had to leave S2S and use an additional result of Büchi \cite{Buchi1973}, stating that the MSO theory of 
$(\omega_1,<)$ is decidable.

In \S\ref{sec:Rabin-Borel-determinacy} we use the methods developed for coding Borel sets of arbitrary complexity to give an effective procedure for deciding which player has the winning strategy in the Determinacy Game for a given Borel set. Concluding remarks are given in \S\ref{sec:remarks}.

\section{Well-levelled orders and their basic properties}\label{sec:def}

\begin{notation}\label{pred} (1) For a partial order $(P,\le_P)$ and $p\in P$ we denote 
\[
\le^P_p=\{q\in P:\,q\le_P p\}, <^P_p=\le^P_p \setminus\{p\}.
\]
When no confusion arises, we write $P$ for
$(P,\le_P)$ and $\le$ for $\le_P$.

{\noindent (2)} For $p <_P q$ we write $[p,q]=\{r\in P:\,p\le _P r\le_P q\}$ and $[q,p]=\emptyset$. The notation $[p,p]$ is another way to denote $\{p\}$. For $p, q\in P$ we write $p\sim q$ if:
\begin{itemize}
\item $p$ and $q$ are comparable,
\item $[p, q]\cup [q,p]$ is a well-order and
\item for any $r\notin [p, q]\cup [q,p]$, we have $r\le_P p$ iff $r\le_P q$.
\end{itemize}
\end{notation}

\begin{observation}\label{obs:equiv} The relation $\sim$ is an equivalence relation on chains.
\end{observation}

\begin{proof} The relation is clearly reflexive and symmetric. Suppose that $p\sim q$ and $q\sim r$ are three distinct mutually comparable elements of
a partial order $P$.

If $p\le q$ and $q\le r$, then clearly $p\le r$ and $[p,r]=[p,q]\cup [q,r]$ is a well-order. If $s\le r$ and $s\notin [p,r]$, then
$s\notin [q,r]$ and so $s\le q$ by the definition of $q \sim r$. But $s\notin [p,q]$ so $s\le p$ by the definition of $p\sim q$.

Suppose that $p\le q$ and $r\le q$. If $r\in [p,q]$ then $p\le r$ and $[p,r]$ is a well-order as a suborder of $[p,q]$. If $s\notin [p,r]$ and $s\le p$, then $s\le r$ since $p\le r$. If $s\notin [p,r]$ and $s\le r$, then $s\notin [p,q]=[p,r]\cup [r,q]$ and $s\le q$. Hence by $p\sim q$ we have $s\le p$.

If $r\notin [p,q]$, since $[q,p]=\emptyset$, we have $r\notin [p, q]\cup [q,p]$. Since 
$r\le q$, it must be that $r\le p$ by the definition of $p \sim q$. In this case we have $[r,p]\subseteq [r,q]$, so $[r,p]$ must be a well-order. If $s\notin [r,p]$
and $s\le p$ then $s\le q$. It cannot be that $s\in [r,q]$ since then $s \in [r,p]$. Therefore $s\notin [p,q]$ and hence $s\le r$, as $s\le q$, 
$s\notin [r,p]$ and $q\sim r$.

The other two cases are similar.
\begin{comment} Suppose $\q\le p$ and $q\le r$. If $p\le r$ then $[p,r]\subseteq [q,r]$ is a well-order. If $s\notin [p,r]$ and $s\le p$, then clearly $s\le r$.
If $s\le r$, if $s\notin [q,r]$, then $s\le q$ since $q\sim r$. It follows that $s\le p$. If $s\in [q,r]$, since $s\notin [p,r]$, then $s\in [q,p]$ and hence $s\le p$.

Otherwise $r\le p$ since we have assumed that the elements $p,q,r$ come from a chain. Then by transitivity we have $p=r$, which is a contradiction with the assumption that the three elements are distinct.

The final case is $q\le p$ and $r\le q$. Hence $r\le p$ and $[r,p]$ is well-ordered as a subset of $[r,q]$. If $s\le r$ then $s\le p$ by transitivity. If 
$s \notin [r,p]$ and $s\le p$, then $s\le q$ by transitivity and $s\notin [r,q]=[r,p]\cup [p,q]$. Since we have $r\sim q$, we obtain $s\le r$.

\end{comment}
\end{proof}

\begin{definition}\label{main-def} (1) Let $(P,\le_P)$ be a partial order. By induction on an ordinal $\alpha$ we define what it means for $P$ to have the
{\em well-level degree} $\le\alpha$, which we denote as $\wl(P)\le\alpha$. We let $\wl(\emptyset)=-1$.

\underline{$\alpha=\beta+1$, possibly $\beta=-1$.} For every $p\in P$, there is $q\sim p$ such that $q\le_P p$ and
$\wl(< ^P_q)\le\beta$. 

 \underline{$\alpha$ a limit $>0$.} For every $p\in P$, $\wl(< ^P_p)\le\beta$ for some $\beta<\alpha$.

\smallskip

{\noindent (2)} If there is an ordinal such that $\wl(P)\le\alpha$, we say that $P$ is {\em well-levelled}. If $P$ is not well-levelled we write $\wl(P)=\infty$.
 If $P$ is well-levelled we say $\wl(P)=\alpha$ iff $\alpha$ is the least ordinal such that $\wl(P)\le\alpha$.
\end{definition}

\begin{example}\label{ex-basic} (1) Linear orders $P$ with $\wl(P)=0$ are exactly the well-orders.

This is because if $q$ is such that $\wl(< ^P_q)= -1$, it means that $\wl(< ^P_q)=\emptyset$, so $q$ is a minimal element of $P$. Therefore the assumption
$\wl(P)=0$ gives that every $p\in P$
satisfies $p\sim q$ for some $q$ which is $P$-minimal. Since $P$ is linear, there is only one minimal $q$ and moreover every $p_0, p_1\in P$ are
comparable. Hence $P$ is the well-ordered chain $\bigcup_{p\in P} [q,p]$.

Obviously, if $P$ is a well-order then $P$ is well-levelled and $\wl(P)=0$, by definition.

{\noindent (2)} $\wl(P)=0$ iff $P$ is a disjoint union of well-orders with no relation between them.

In the backward direction, suppose $P=\bigcup_\alpha C_\alpha$ and each $C_\alpha$ is a well-order, while there is no relation between the elements of
$C_\alpha$ and $C_\beta$ for $\alpha\neq \beta$. Then each $C_\alpha$ consists of $\sim$-equivalent elements and $\wl(C_\alpha)=0$, so $\wl(P)=0$.

In the forward direction, let us write $P$ as the union of $\sim$-equivalence classes, which we shall enumerate as $\{C_\alpha:\alpha<\kappa\}$. To see that this is possible, apply the argument of (1) to conclude from $\wl(P)=0$ that every element of $P$ is above a unique minimal element and $\sim$-equivalent to it.
Hence we can let $C_\alpha$ be the set of elements equivalent to a given minimal element, and
then each $C_\alpha$ is a linear order, since $p\sim q$ implies that $p$ and $q$ are comparable. For $\alpha<\kappa$ let $q_\alpha$ be 
the unique minimal element of $C_\alpha$. Then we conclude as in (1) that $C_\alpha$ is the well-order $C_\alpha=
\bigcup_{p\in C_\alpha} [q_\alpha, p]$. We shall show that there is no relation between the elements of $C_\alpha$ and $C_\beta$
for $\alpha\neq \beta$.

Suppose for a contradiction that $p\in C_\alpha$, $q\in C_\beta$ and $p\le q$. Since $p \notin [q_\beta, q]$, we have $p\le q_\beta$, but this is a contradiction
with $q_\beta$ being $P$-minimal and $p\neq q_\beta$.

{\noindent (3)} Let $p_l (l<3)$ be such that $p_0< p_1, p_2$ and $p_1$ and $p_2$ are incomparable. Let $P=\{p_l:\;l<3\}$. Then $\wl(P)=1$.

{\noindent (4)} A more general example than that of (3) is a tree, a particular instance of which is the order in (3). We claim that if $T$ is any tree which is not a linear order, then $\wl(T)=1$.

Namely, if $T$ is rooted with the root $r$, then we have $\wl(\{r\})=0$. Every $p\in P$ satisfies 
$p\sim r$, so $\wl(T)=1$. If $T$ is not a rooted tree, it is a disjoint union of rooted trees with no relations between them. Since each of them is well-levelled and has $\wl=1$, we have $\wl(T)=1$.

{\noindent (5)} Let $p_l (l<3)$ be such that $p_0,  p_1 < p_2$ and $p_0$ and $p_1$ are incomparable. Let $P=\{p_l:\;l<3\}$. Then $P$ is well-levelled and 
$\wl(P)=1$.

\end{example}

\begin{proposition}\label{prop:basic-properties}  {\noindent (1)} If $P$ is a well-levelled order and $Q\subseteq P$ is a sub-order which is downward closed, then $Q$ is
well-levelled and $\wl(Q)\le \wl(P)$.

 {\noindent (2)} Every well-levelled order is well-founded.

\end{proposition}

\begin{proof} (1)  Let $P$ be a well-levelled order with $\wl(P)=\alpha$.  Notice that for every $q\in Q$
we have $<^q_Q=  <^q_P$, since $Q$ is downward closed. Similarly, for $p,q\in Q$, we have $[p,q]^Q=[p,q]^P\cap Q$. In particular, we have that 
$(p\sim_P q) \implies (p\sim_Q q)$. The proof now proceeds by induction on $\alpha$.

Suppose that \underline{$\alpha=0$}. By Example \ref{ex-basic}(2) we have that $P$ is a disjoint union of well-ordered chains with no relation between them, so then clearly so is
$Q$ and hence $\wl(Q)=0$.

Suppose \underline{$\alpha=\beta+1$} for some $\beta\ge 0$. Let $p\in Q$, so by hypothesis about $P$, there is $q\in P$ with $q\sim_P p$, $q\le_P p$ and 
$\wl(<^P_q)\le \beta$. By the hypothesis that $Q$ is downward closed, we have $q\in Q$ for any such $q$. We also have $q\sim_Q p$ and $<^Q_q= <^P_q$, hence we are done by the induction hypothesis.

If \underline{$\alpha>0$ is a limit}, then the conclusion follows directly by the downward closure of $Q$.

{\noindent (2)} Let $P$ be a well-levelled order with $\wl(P)=\alpha$. The proof is by induction on $\alpha$.

The case \underline{$\alpha=0$} follows by Example \ref{ex-basic} (2). 

In the case \underline{$\alpha=\beta+1$} for $\beta\ge 0$, suppose for a contradiction that there is an infinite decreasing sequence $\langle p_n:\,n<\omega\rangle$ in $P$. Then $\langle p_n:\,1\le n<\omega\rangle$ is an infinite decreasing sequence in $<^P_{p_0}$, while we know by definition of
$\wl(P)\le\beta+1$ that there is some $q\sim p_0$ with $q\le p_0$ such that $\wl(<^P_q)\le \beta$. In particular, the intersection of $\langle p_n:\,1\le n<\omega\rangle$
with $[q,p_0]$ is finite, since $[q,p_0]$ is a well-order. By the definition of $\sim$, for every
$n$ such that $p_n$ is not in $[q,p_0]$ we have that $p_n\le q$. Hence the intersection of $\langle p_n:\,1\le n<\omega\rangle$
with $<^P_q$ is infinite, a contradiction with the induction hypothesis.

The proof in the case \underline{$\alpha>0$ a limit ordinal} is similar to the successor case.
\end{proof}

To give the reader a bit of the intuition, the following Observation \ref{ex-basic2} gives some examples and easy remarks.
The basic intuition is that a well-order is well-levelled of degree 0, a tree is a partial order in which every element has the set of predecessors which is a well-order of level 0, so $T$ is a well-levelled order of degree 1. And so on, the examples are numerous and we have basic operations that give us more well-levelled orders out of the given ones.

\begin{observation}\label{ex-basic2}
{\noindent (1)} Let $T$ be a tree and $p\notin T$. Define $P=T\cup \{p\}$ and let $\le_P\rest T=\le_T$, while $t\le_P p$ for every $t\in T$. Then $\wl(P)=2$.

{\noindent (2)} Let $T$ be a tree and $p\neq q$ two objects not elements of $T$. Define $P=T\cup \{p,q\}$ and let $\le_P\rest T=\le_T$, while $t\le_P p, q$ for every $t\in T$
and $p\le_Pq$. Then $\wl(P)=2$. 

{\noindent (3)} Let $T_0, T_1$ be disjoint rooted trees and $P$ the lexicographic sum $T_0+T_1$, that is the order in which all elements of $T_0$ are below any element of $T_1$ and the orders in $T_0$ and $T_1$ are induced by the original orders. Then $\wl(T_0+T_1)=2$, which can be seen by noticing that
for every $p\in T_1$, we have $p\sim r$ where $r$ is the root of $T_1$. Then we reduce the consideration to (1).

In fact, the disjoint union of any two well-levelled orders with no relations in between is well-levelled. Indeed, if $P$ and $Q$ are such orders then for every $r\in R=P\cup Q$ we have that $<^R_r$ is either $<^P_r$ or $<^Q_r$, both of which are well-levelled by the assumption. A similar conclusion holds for any
number of disjoint parts of a union of well-levelled orders.

{\noindent (4)} For $i<n$ let $T_i$ be a rooted tree and consider the lexicographic sum $T_0+ T_1+\ldots + T_{n-1}$. Then 
$\wl(T_0+ T_1+\ldots + T_{n-1})=n$.

{\noindent (5)} The following is a main operation for the purposes of this paper. It will be generalised in Definition \ref{def:powers}. The author used this construction in \cite{zbMATH01310540} to solve a question in combinatorics (posed by R. Ahlswede,  Peter Erd\H{o}s and N. Graham) .

\begin{notation}\label{not:Tsquare}
Let $T$ be a rooted tree of height $\omega$ and replace every $t\in T$ by a copy $T_t$ of $T$, rooted at $t$. Let $P=\bigcup_{t\in T} T_t$ and for $p,q\in P$ declare
$p\le_P q$ if either $p, q$ are in the same copy $T_t$ of $T$ and $p\le_T q$, or $p\in T_s, q\in T_t$ and $s<_T t$. Denote this operation by $T\cdot T$
and the resulting order by $T^2$. 
\end{notation}

We claim that 
$\wl(T^2)= \omega$. 

To see this, let $p\in P$ and let $s\in T$ be such that $p\in T_s$. Then $p\sim r^s$, where $r^s$ is the root of $T_s$. Let $n={\mbox ht}_T(s)$ We have that $<_{r_s}^P$ is isomorphic to a lexicographic sum of $n$ trees, and hence by (4), $\wl(<_{r_s}^P)=n$. Therefore $\wl(T^2)= \omega$. 

The elements of $T^2$ can be given coordinates as pairs $(s,t)$ where $t$ says that we are inside of the tree $T_t$ and $s$ that we have taken $s\in T_t$.
(The reader might find the order of taking $t$ as the second coordinate somewhat unnatural, but this is justified by our needs in Definition \ref{def:powers}.)

We can similarly define the operation $P\cdot T$ where $P$ is a well-levelled partial order and $T$ is a rooted tree. The next two examples will give
the instances of this operation that are most relevant to us.

{\noindent (6)}  Let $T$ be a rooted tree and we define by induction on $n\ge 1$ what is $T^n$. Let $\Lambda$ denote the root of $T$.

$T^1=T$ and for $n\ge 1$ we have $T^{n+1}=T^n\cdot T$, that is, every node in $T$ is replaced by a copy of $T^n$.

Notice that for $n\ge 1$ the elements of $T^n$ can be given the coordinates in the $n$-fold Cartesian product of $T$ with itself. The projections in this
notation are defined naturally, where we use the notation $x \rest A$ for the projection of an element $x$ on a subset $A$ of the coordinates and $x(k)$ for the value of $x$ at the
coordinate $k$. We now define the order $\le_{T^n}$ in $T^n$, for $n\ge 1$.
The definition is by induction on $n$, as follows:
\begin{itemize}
\item $\le_{T^1}=\le_T$,
\item for $x,y\in T^{n+1}$ we let $x\le_{T^{n+1}} y$ if $x(n)=y(n)$ and $x \rest {[1,n)}\le_{T^n} y \rest {[1,n)}$, or $x(n)<_T y(n)$.
\end{itemize}

Let $T^{<\omega}=\bigcup_{1\le n<\omega} T^n$ ordered by the induced order. Within $T^{<\omega}$ we have the natural projections maps 
$p^n_k:\, T^{n}\to T^k$ for $k\le n$, where $p^n_k(x)=x\rest [1,k)$.

We define $T^{\omega}$ as the set of all $\omega$-sequences $x$ where there is some $n$ such that $x\rest [1, n)\in T^n$ and for all $m\ge n$ we have
$x(m)=\varLambda$. For $x\in T^\omega$ let $n(x)$ be the first $n$ such that $x(n)=\lambda$ for all $m\ge n$. We define the order on $T^\omega$ by letting
$x\le y$ iff $n(x)=n(y)$ and $x\rest [1, n(x)) \le_{T^{n(x)}} y \rest [1, n(y))$. This definition allows us to consider elements of $T^\omega$ as sequences of length $\omega$ of elements of $T$, which take value $\Lambda$ on a co-finite set of coordinates.

{\noindent (7)} 
By convention, we set $T^0=\emptyset$ and we do not define $T^{<0}$. For a successor ordinal $\alpha=\beta+1$ we define $T^\alpha=T^\beta\cdot T$
and for $x,y\in T^{\alpha}$ we let $x\le_{\alpha} y$ if $x(\beta)=y(\beta)$ and $x\rest [1,\beta)\le_{T^\beta} y\rest [1,\beta)$, or $x(\beta)<_T y(\beta)$. \end{observation}

By analogy with examples (6) and (7) and with the same notation, we can define $(T^\alpha, \le_{T^\alpha})$ for any ordinal $\alpha>0$, as follows.
 
\begin{definition}\label{def:powers} Let $T$ be a rooted tree with the root $\Lambda$. By induction on $1\le\alpha$ we define what is $T^{\alpha}$, which will be a subset of the set of functions ${}^{\alpha} T$. We also define the associated partial order $\le_{T^\alpha}$.
The definition is by induction on $\alpha$.

For \underline{$\alpha=1$} we have that $T^1=T$ and $\le_{T^1}=\le_T$. 

For \underline{$\alpha=\beta+1$} we have that $T^{\beta+1}$ is a copy of $T$ in which every element of $T$ is replaced by a copy of $T^\beta$.
Formally,
\[
T^{\beta+1}=\{(s,t):\,s\in T\wedge t\in T^\beta\},
\]
where we have identified any function $t\in {}^{\alpha} T$ with the pair $(t(\beta), t\rest [1,\beta))$. The order $<_{T^{\beta+1}}$
is given by $(s,t)<_{T^{\beta+1}}(s',t')$ if $s=s'$ and $t<_{T^\beta}t'$ or $s<_T s'$.

For \underline{$\alpha>0$ a limit} we let
\[
T^\alpha=\{t \in {}^{\alpha} T:\,(\forall \beta<\alpha)[ t\rest \beta\in T^\beta]\wedge (\exists \beta<\alpha)[ t\rest [\beta, \alpha)\equiv \Lambda]\}.
\]
This completes the definition of $T^\alpha$, and we now define the partial order $\le_{T^\alpha}$. For $s, t\in T^\alpha$ we let $s\le_{T^\alpha} t$ iff for the first $\beta<\alpha$ such that 
$s\rest [\beta, \alpha)=t\rest [\beta, \alpha)\equiv\Lambda$ we have $s\rest\beta\le_{T^\beta} t\rest\beta$.
\end{definition}

\begin{remark}\label{rem:moregenerality} Definition \ref{def:powers} is in fact nothing else but the definition of an inverse limit of the products of $T$. It can clearly be extended to the more general case when not all trees in the limit are the same or when the factors are not trees but just any well-levelled orders with
the least element. Since in this paper we shall only be interested in the case that $T$ is the Cantor tree $2^{<\omega}$ and $\alpha<\omega_1$, we 
do not develop the definitions in the unnecessary generality.
\end{remark}

\begin{notation}\label{not:restrictions} From this point on, we use the notation $T$ for the Cantor tree $2^{<\omega}$ whose root is denoted by $\Lambda$.

If $A\subseteq T^\alpha$ and $\beta<\alpha$, we denote $A\rest[1,\beta)=\{x\rest[1,\beta):\,x\in A\}$.
\end{notation}

\begin{lemma}\label{lem:Cantorproducts} The order $T^\alpha$ for $1\le \alpha<\omega_1$ is well-levelled and countable. 
\end{lemma}

\begin{proof} The proof is by induction on $\alpha$. 

The case \underline{$\alpha=1$} is clear, since $T$ is the Cantor tree.

Suppose that \underline{$\alpha=\beta+1$} where $\beta<\omega_1$. By the induction hypothesis we have that $T^\beta$ is countable, so then 
$¨T^\alpha=T^\beta\cdot T$ is countable as well. 

Given $y\in T^\alpha$, we consider $<^{T^\alpha}_y$. It is the disjoint union 
\[
\{x\in T^\alpha:\,x(\beta)=y(\beta)\wedge x\rest [1,\beta)<_{T^\beta} y\rest[1,\beta)\}\cup \{x\in T^\alpha:\,x(\beta)<_T y(\beta)\}.
\]
By Observation \ref{ex-basic2}(3), it suffices to show that each of the parts of the union is well-levelled. The first part is isomorphic to 
$<^{T^\beta}_{y\rest[1,\beta)}$, which is well-levelled since $T^\beta$ is. The second part is isomorphic to $T^\beta\cdot <^T_{y(\beta)}$. We have that
$<^T_{y(\beta)}$ is a finite linear order, say of the finite size $n$. Therefore $T^\beta\cdot <^T_{y(\beta)}$ is isomorphic to $T^\beta\cdot n$. The conclusion
follows from Observation \ref{ex-basic2}(5) where the ``tree'' in question is just the linear order with $n$ elements and $P$ is $T^\beta$.

Suppose that \underline{$0<\alpha <\omega_1$ is limit}. Every element of $T^\alpha$ is eventually $\Lambda$, so there is an injection of $T^\alpha$
into $\bigcup_{\beta<\alpha}T^\beta$, which is countable by the inductive hypothesis. Therefore $T^\alpha$ is countable. If $y\in T^\alpha$, let 
$\beta<\alpha$ be the first such that $y(\gamma)=\Lambda$ for every $\gamma\in[\beta,\alpha)$. Hence $<^{T^\alpha}_y$ is isomorphic to 
$<^{T^\beta}_{y\rest[1,\beta)}$, which is well-levelled by the inductive assumption.
\end{proof}

\begin{definition}\label{def:Boreltres} Let $\alpha<\omega_1$ and let $T$ be the Cantor tree $2^{<\omega}$. 
%We shall sometimes abuse the notation by denoting the order relation in $T^{1+\alpha}$ by $\subseteq$.

A {\em path} in $T^\alpha$ is a linearly ordered subset $\pi$ such that for no $\beta\le\alpha$ does $\pi\rest\beta$ have an upper bound. Note that by the fact that $T^\alpha$ is well-founded, every path is
order isomorphic to an ordinal. An {\em $\omega$-path} is a path that is order isomorphic to $\omega$. A {\em branch} is a path that is downward closed. (For 
$\alpha=1$ this coincides with the maximal paths, but not for $\alpha$ is general.)
The set of all branches in $T^\alpha$ is denoted by $[T^\alpha]$.

In the case of an $\omega$-path $\rho$, we denote by $\rho\rest k$ the chain consisting of the first $k$ elements of $\rho$, for $k<\omega$.
\end{definition}

\subsection{The Polish Topology of $[T^{1+\alpha}]$}\label{sec:metric} This subsection is not used in the proof of Theorem \ref{th:main}. We feel that it is useful for the intuition of in what sense the coding of Borel sets developed in \S\ref{sec:realcoding} resembles the usual coding of closed sets in the Cantor set $2^\omega$ by the branches of subtrees of $T$. We code closed and open sets in the Polish space $2^\omega=[T]$ as sets of branches in the tree $T$. Similarly, we
shall code $\Sigma_\alpha\cup \Pi_\alpha$-sets in $2^\omega$ as sets of branches in the well-levelled order $T^{1+\alpha}$ (see \S\ref{sec:realcoding}). The set $[T^{1+\alpha}]$ of all branches of $T^{1+\alpha}$ is itself a Polish space. Moving from the original Polish space to different Polish spaces for the codes is one of the main differences between what is done here and what was done previously for the closed sets. Rabin's idea for dealing with the $F_\sigma$-sets can a posteriori also be understood as a step in this direction. Let us therefore describe these Polish spaces.

For countable ordinals $\alpha$, each $[T^{1+\alpha}]$ has a naturally defined complete 
metric, so the resulting topology, which we shall denote by $\tau_\alpha$ makes $[T^{1+\alpha}]$ into a Polish space. 
To start the proof of this, recall (\cite[pg. 7]{Kechris}) that the Cantor space $[T]$
has a complete metric given by $d_0=d$ where,
\[
d(f,g)=\begin{cases}
2^{(-\min\{n :\,f(n)\neq g(n)\}+1)} &\text{if }f\neq g \\
0 &\text{if } f=g.
\end{cases}
\]
This metric is bounded by 1. The following classical theorem shows that this metric can be used to make any countable product $\prod_{1+\alpha} [T]$
into a Polish space. It is proved in \cite[Th. 4.3.12]{E}.

\begin{theorem}[{\bf The Polish Prewar School}]\label{th:product-metric} Suppose that $(X_n, \psi_n)_{n<\omega}$ is a sequence of non-empty complete metric spaces with metrics bounded by
$1$. Then the product $\prod_{n<\omega} X_n$ admits a complete metric given by
\[
\rho(f,g)=\sum_{n<\omega} 2^{-(n+1)} \psi_n(f(n), g(n)).
\]
\end{theorem}

There is more than one way to define a compatible metric for a given Polish topology. For example \cite[pg. 3]{Kechris} gives a different formula for the
metric of a countable product of completely metrisable spaces. We like the one given in the statement of Theorem \ref{th:product-metric} because when applied to the countable product of the 0-1 metrics on $2=\{0,1\}$, it gives the natural metric on $2^\omega$.

Theorem \ref{th:product-metric} can be adapted to show the following Lemma \ref{lemma:Yair}.\footnote{We thank Yair Hayut who during our talk on this paper at the Set Theory Seminar at the Hebrew University of Jerusalem on 14/01/2026 noticed a mistake in the previous version of the proof of Lemma \ref{lemma:Yair}. The mistake was due to the definition of a branch which did not fit what we have considered a branch in the rest of the paper.}

\begin{lemma}\label{lemma:Yair} For every $\alpha<\omega_1$, the space $[T^{1+\alpha}]$ admits a complete compatible metric $d_\alpha$ bounded by 1. 
\end{lemma}

\begin{proof} The proof is by induction on $\alpha$. 

For \underline{$\alpha>0$}, this is the classical fact that the Cantor set $[T]$ has a complete metric which generates the topology, while we have 
$[T^1]=[T]$.

For \underline{$\alpha=\beta+1$}, first note that if $\rho$ is a branch in $T^{1+\alpha}$, then $\rho\rest [1+\beta)$ is a branch in
$T^{1+\alpha}$, by the definition of what a branch is. 

We define $d_\alpha(f,g)=\frac{1}{2} \cdot [d_\beta(f\rest [1,\beta), g\rest [1,\beta)) + d_1(f(\beta), g(\beta))]$. This function is clearly bounded by 1, non-negative and symmetric, and $d_\alpha(f,g)=0$ iff $f=g$. The triangle inequality follows by simple substitution and the use of the triangle inequalities for
$d_\beta$ and $d_1$.

To see that the metric is complete, note that if $\langle f_n:\,n<\omega\rangle$ is a Cauchy sequence in $d_\alpha$, then 
$\langle f_n\rest [1, \beta):\,n<\omega\rangle$ is a Cauchy sequence in $d_{\beta}$, hence it converges to some $f\in [T^{1+\beta}]$. Similarly,
$\langle f_n\rest (\beta):\,n<\omega\rangle$ is a Cauchy sequence in $d_1$, so it converges to some $g\in [T]$. Then 
$\langle f_n:\,n<\omega\rangle$ converges to the branch $\{t\!\frown \! s:\,t\in f \wedge s\in g\}$ in $[T^{1+\alpha}]$. 

For \underline{$\alpha>0$}, let $\langle \beta_n:\,n<\omega\rangle$ be an increasing cofinal sequence in $\alpha$. Let 
$f,g \in  [T^{\alpha}]$.
Define $d_\alpha(f,g)$ using the formula similar to that of Theorem \ref{th:product-metric} with $\rho_n=d_{\beta_n}$, namely
\[
d_\alpha(f,g)=\sum_{n<\omega} 2^{-(n+1)} d_{\beta_n}(f\rest [1,\beta_n), g\rest [1,\beta_n)),
\]
but note that $f$ and $g$ here are not the elements of any product, but branches in $[T^{\alpha}]$. The definition is justified since for every
$f\in [T^{1+\alpha}]$ we have $f\rest [1+\beta_n)$ is a branch in $T^{1+\beta_n}$.

The proof that this is a metric bounded by 1 is similar as in the successor case. To show that it is complete, suppose that 
$\langle f_n:\,n<\omega\rangle$ is a Cauchy sequence in $d_\alpha$. It follows that for every $\beta_n$, the sequence 
$\langle f_n\rest [1, \beta):\,n<\omega\rangle$ is a Cauchy sequence in $d_{\beta_n}$ and hence converges to some $\pi_n \in [T^{1+\beta_n}]$.
Moreover, it follows by the definitions taken for $d_{\beta_n}$'s, that for $n\le m$ we have
$\pi_n=\{t\rest [1+\beta_n):\,t \in \pi_m\}$. Let 
\[
\pi=\{t\in T^{\alpha}:\,(\forall n) t\rest [1+\beta_n) \in \pi_n \}.
\]
We claim that $\pi$ is an element of  $[T^{\alpha}]$ and the limit of $\{\pi_n:\,n<\omega\}$. First let us show that it consists of pairwise compatible elements in $T^{\alpha}$, so suppose that $t,r\in \pi$ and let $k$ be the smallest such that $t\rest [\beta_k, \alpha)=r \rest [\beta_k, \alpha)\equiv \Lambda$. 
By the definition, $t\le_{T^{\alpha}}r$ iff $t\rest [1+\beta_k)\le_{T^{1+\beta_k}} r\rest [1+\beta_k)$ and vice versa. One of these possibilities must occur since 
$\pi_k$ is a branch in $T^{1+\beta_k}$. We similarly see that $\pi$ is downward closed. Let us calculate $d_\alpha(\pi, \pi_n)$ for an arbitrary $n<\omega$. By
the definition, this is $\le \sum_{n\le m<\omega} 2^{-(m+1)}$, a quantity that converges to 0 as $n$ increases. Hence $\pi$ is the limit, as required to show.
\end{proof}

\section{Coding of the Borel Sets}\label{sec:realcoding} Recall that $T$ stands for the Cantor tree $2^{<\omega}$. To every Borel set $H \subseteq 2^\omega$ we shall associate the code 
$c(H)$ which will be a set of branches of a structure contained in $\bigcup_{0< \alpha<\omega_1}T^\alpha$. Recall that the Borel sets are the members of
the hierarchy of sets defined by the increasing sets of levels $\Sigma_\alpha\cup \Pi_\alpha$ for $\alpha<\omega_1$ (see \cite[\S11.B]{Kechris}). In Definition \ref{def:Borelcodes} we use the
induction on $\alpha<\omega_1$ to describe the coding of the sets $H$ in $\Sigma_{\alpha}\cup \Pi_\alpha$, which will be done using sets of branches
in $T^{1+\alpha}$. 

To start the definition, we recall the well-known coding of the closed subsets $H$ of the Cantor set $2^\omega$. We identify the product $2^\omega$
of $\omega$ many copies of $2=\{0,1\}$  in the product topology with a topological space obtained on the set ${}^\omega 2$ of $\omega$-sequences
of 0s and 1s, with a basis of the topology
given by the sets $[s]=\{f\in {}^\omega 2:\,s\subseteq f\}$ for $s$ a finite partial function from $\omega$ to 2 (equivalently, we can restrict our attention to the $s$ whose domain is some natural number as a sub-basis). It is usual to denote both of these objects by
$2^\omega$ and call them both the Cantor set. 

A branch $\rho$ in the Cantor tree $T=2^\omega$ is identified with the element $\bigcup \rho$ of the Cantor space $2^\omega$, and it is usual to denote both of these by $\rho$ and talk about them interchangeably. Then it is well known and easily verified on the basis of the definitions given so far, that any non-empty closed set $H\subseteq 2^\omega$ is exactly (identified with) the set $[S]$ of branches of a subtree $S$ of $T$ of height $\omega$, and vice versa, any such set of branches forms a closed set. See \cite[\S2.B]{Kechris} for a detailed exposition. We use the notation $T_H$ for a tree $S$ obtained from the closed set
$H$ in this manner. In other words, $T_H=\{s\in T:\,[s]\cap H\neq \emptyset\}$.
\begin{comment} Let $H$ be a non-empty closed subset of $2^\omega$ and let $S=\{f\rest n:\,f\in H \wedge n<\omega\}$. Then $S\subseteq T$ is downward closed by definition, and it has height $\omega$ since there is at least one element in $H$. If $g\notin H$ then there is the least $n$ such that $g\rest n\notin S$. It follows that the basic open set $[g\rest n]$ is disjoint from $[S]$. Hence the complement of $[S]$ is open and so $[S]$ is closed. We clearly have
$H\subseteq [S]$. 

We still need to show that $[S]\subseteq H$, so suppose that $g\in 2^\omega$ is such that for all $n<\omega$ we have $g\rest n\in S$. Therefore,
for each $n<\omega$ there is $f_n\in H$ such that $f_n\rest n=g\rest n$. It follows that $d(g,f_n)\le 1/n$ in the usual metric of the Cantor space. In particular,
$g$ is the limit of a sequence $\langel f_n:\,n<\omega\rangle$ of elements of $H$. Hence $g\in H$, since $H$ is closed.
\end{comment}

We then may say that open sets are the sets of the form $2^\omega\setminus [T_H]$ for some $H$ closed. We therefore obtain that both open and closed sets are identified with, which we will name ``coded by'', exactly themselves.

\begin{notation}\label{note:projection-base} (1) For $\alpha<\omega_1$ and $x\in T^{1+\alpha+1}$, we denote by:
\begin{itemize}
\item $p(x)$ ($p$ stands for {\em projection}) the element of $T^{1+\alpha}$ such that $x= p(x)\frown y$ for some $y\in T$ and 
\item $b(x)$ ($b$ stands for {\em base}) the $y\in T$ such that $x=p(x)\frown y$.
\end{itemize}

{\noindent (2)} For a branch $\pi$ of $T^{1+\alpha+1}$ and $n<\omega$, we denote 
\[
\pi^n=\{ p(x): x\in \pi \wedge {\rm ht}(b(x))=n\}.
\]
%
\begin{comment}
{\noindent (3)} Suppose that $\rho$ is an $\omega$-path in $T^{1+\alpha}$ for some $\alpha<\omega_1$ and that $n<\omega$. Note that the projection of
$\rho$ to any $\{\beta\}$ for $\beta<\alpha$ is an $\le\omega$-path in $T$. We denote by $\rho(\beta)$ that path of $T$.
We denote by
$\rho \Rsh n$ the $1+\alpha$ sequence of paths such that for every $\beta<1+\alpha$ we have $(\rho \Rsh n)(\beta)=\rho(\beta)\rest n$.

{\noindent (4)} Suppose that $\beta\le\alpha$ and $\rho$ is a path of  $T^{1+\alpha}$. By the {\em projection} of $\rho$ to $1+\beta$
we mean the set $\{t\rest (1+\beta):\,t\in \rho\}$, which is then a  path of $T^{1+\beta}$. Following Notation \ref{not:restrictions}, we denote this projection by
$\rho\rest (1+\beta)$.

{\noindent (5)} We shall use the expression {\em MSO-definable} to refer to MSO-definitions with respect to the Cantor tree and the 2 successor relations, the S2S theory of \cite[Definition 1.1]{Rabin} (the two orders on $T$, the natural and the lexicographic one, are actually MSO-definable from the two successor relations, see \cite{Rabin}).
\end{comment}
%
\end{notation}

\begin{remark}\label{rem:restrictions} We emphasise that the notation for restriction $\rest$ is used for different objects with different results, depending on the
domain of the function to restrain. So, in the case of $\rho$ a path in $T$, we consider it as a function from a subset of $\omega$ to 2 and hence
it makes sense to talk about $\rho\rest n$. 

In the case of $t\in T^{1+\alpha}$ we consider $t$ as a function from
$1+\alpha$ to $T$ and hence for $\beta\le\alpha$ the notation $t\rest (1+\beta)$ makes sense and denotes the function with domain $1+\beta$ whose value at $\gamma\le 1+\beta$ is $t(\gamma)$. In the case of $\rho$ a path in $T^{1+\alpha}$, we consider it as a set of functions with domain $1+\alpha$, 
and the notation $\rho\rest(1+\beta)$ for $\beta<\alpha$ is appropriate as well, as per Notation \ref{not:restrictions}.

Within the following Definition \ref{def:Borelcodes} we shall briefly indicate how the choices done in it can be done in an MSO-definable way. This is not yet needed in this section but is crucial for the main result of the paper Theorem \ref{th:main}. We devote the full section \S\ref{sec:automates} to explaining this properly, but it seems opportune to mention this issue already at this point as a motivation for the long definition. Presenting this definition separately and before the methodology of MSO-definability in $T^{1+\alpha}$ is an editorial choice which in our opinion improves the readability of the paper and emphasises the main idea. See the penultimate paragraph of \S\ref{sec:remarks} for how a different approach, involving a long simultaneous induction between the definition here and the proofs in \S\ref{sec:automates} would have to be done for computer verification purposes.

\end{remark}

\begin{definition}[{\bf Coding Borel sets}]\label{def:Borelcodes} By induction on $\alpha<\omega_1$ we describe the coding of the sets $H$ in the level
$\Sigma_{\alpha}\cup \Pi_\alpha$ of the Borel hierarchy of the Cantor space. Such sets will be coded using sets of branches
in $T^{1+\alpha}$. 

\underline{$\alpha=0$.}  For a closed set $H$  which is either open or closed, we let $c(H)=H$, hence this code is a set of branches in $T=T^{1+0}$. 

\underline{$\alpha+1$.} Given a non-empty $\Sigma_{\alpha+1}$-set $H=\bigcup_{n<\omega} H_n$, where each $H_n$ is a $\Pi_\alpha$-set which we
presume coded by a set $c(H_n)$ of branches in $T^{1+\alpha}$. We shall only proceed if the sequence 
$\langle H_n:\,n<\omega\rangle$ is $\subseteq$-increasing and $H_0\neq \emptyset$. This can be done without loss of generality since any non-empty 
$\Sigma_{\alpha+1}$-set has a representation satisfying that property. Let $B$ be the branch of $T$ consisting of constantly $0$ sequences, including the
empty one $\lambda$. Let $c(H)$ be the set of all branches $\rho^\ast$ of $T^{1+\alpha+1}$ where for every $n<\omega$ we have that $\rho^\ast(\alpha)(n)=(0)_{n}$ and $\{t\rest [1+\alpha): \,t \in \rho^\ast \mbox{ and } b(t)= (0)_{n}\}=\rho$ for some element $\rho$ of $c(H_n)$. 

Given a non-empty $\Pi_{\alpha+1}$-set $H=\bigcap_{n<\omega} H_n$, where each $H_n$ is a $\Sigma_\alpha$-set and given $c(H_n)$ a set of branches of some subset $S_n$ of $T^{1+\alpha}$. We can MSO-define (see \S\ref{sec:automates}) such $S_n$ by taking all initial segments with respect to $T^{1+\alpha}$ of branches in
$c(H_n)$, or we could even take the whole $T^{1+\alpha}$, as we shall only use the branches that are actually in $c(H_n)$ (we prefer calling the set $S_n$ for the ease of visualising the process, but an important point to retain for the future is that $S_n$ can be chosen in an MSO-definable way. Note that from the inductive assumption that $c(H_n)$ is MSO-definable, we can decide for any branch of $T^{1+\alpha}$ if it is in $c(H_n)$ or not). We shall only proceed if the sequence 
$\langle H_n:\,n<\omega\rangle$ is $\subseteq$-decreasing, which we can do without loss generality since any $\Pi_{\alpha+1}$-set has a representation satisfying that property.

We define a substructure $s(H)$ of $T^{1+\alpha+1}$ by taking a copy of $T$ in which every node of height $n$ is replaced by $S_n$.

Recall the notation $\pi^n$ from Notation \ref{note:projection-base} for branches
$\pi$ of $T^{1+\alpha+1}$. Note that each such $\pi^n$ is a branch of $T^{1+\alpha}$. Then
we define $c(H)$ as the set of branches $\pi$ of $s(H)$ such that for every $n$, the projection ${\pi}^{n}$ is an element of $c(H_n)$.

\underline{$\alpha$ limit $>0$.} We code non-empty $\Sigma_\alpha$-sets using the branches of $T^\alpha$. To do so we start with a fixed sequence
$\langle \beta_n:\,n<\omega\rangle$ with $\sup_{n<\omega}\beta_n=\alpha$. Again, for purposes of definability, note that such a sequence can be chosen in an MSO-definable way, since in S2S we have access to all countable ordinals and $\omega$-sequences of them. We can make an enumeration of such sequences, put an order between them and choose the first one that is convenient. See \S\ref{sec:automates} on MSO-definability.

Suppose we are given a non-empty $\Sigma_\alpha$-set $H$. Knowing that the families $\Pi_\beta$-increase with $\beta$ and that our representations of
the form $\langle H_n:\,n<\omega\rangle$ are allowed to have repetitions, we can assume that
$H=\bigcup_{n<\omega} H_n$ for some sequence $\langle H_n:\,n<\omega\rangle$ of sets such that each $H_n$ is a $\Pi_{\beta_n}$-set. Therefore, $c(H_n)$ is a set of branches of some MSO-definable subset $S_n$ of $T^{1+{\beta_n}}$. Again, we shall
assume that $H_n$s are $\subseteq$-increasing and $H_0\neq \emptyset$.

We let $c(H)$ be the set of all branches $\pi$ of $T^{\alpha}$ such that for some
$n$, the branch $\pi\rest (1+\beta_n)$ is an element of $c(H_n)$, and for every $t\in \pi$ we have $t(\beta_n)\in \{\Lambda, \langle 0\rangle\}$ and $t\rest (\beta_n,\alpha)\equiv\Lambda$. 

We code $\Pi_\alpha$-sets similarly. Given a non-empty $\Pi_\alpha$-set $H=\bigcap_{n<\omega} H_n$,
where $H_n$ is a $\Sigma_{\beta_n}$-set for some MSO-definable sequence $\langle \beta_n:\,n<\omega\rangle$ with $\sup_{n<\omega}\beta_n=\alpha$. Therefore, $c(H_n)$ is a set of branches of some MSO-definable subset $S_n$ of $T^{1+{\beta_n}}$. We
assume that $H_n$s are $\subseteq$-decreasing and we let $c(H)$ consist of all branches $\pi$ of $T^\alpha$ such that 
for some
$n$, the branch $\pi\rest (1+\beta_n)$ is an element of $c(H_n)$, and for every $t\in \pi$ we have $t(\beta_n)\in \{\Lambda, \langle 1\rangle\}$ and $t\rest (\beta_n,\alpha)\equiv\Lambda$.
\end{definition}

\begin{example}\label{ex:notunique} The set $c(H)$ is not uniquely determined by $H$. For example, we can have a $G_\delta$-set
$H$ which is represented by two different sequences of closed sets $\langle F_n:\,n<\omega\rangle$ and $\langle F'_n:\,n<\omega\rangle$ 
so that $H=\bigcup_{n<\omega} F_n=\bigcup_{n<\omega} F'_n$. Then the way we define $s(H)$ depends on which of these two representations we have chosen. We reduce the number of codes we can make for a given set, because we consider only the MSO-definable ones, but we have still let the 
choice of the representations as $\sigma$-combinations of basic open sets rather free. The importance of doing this is for Lemma \ref{le:intersections}.
The proof of Theorem \ref{th:main} will use Lemma \ref{le:intersections}, but we shall have to explain why the freeness of the representation does not 
harm the MSO-definability of the whole procedure. This will be done in Lemma \ref{le:representations}.
When speaking of the codes $c(H)$ for a Borel set $H$ we shall tacitly assume that we have a fixed representation of $H$ as a $\sigma$-Boolean
combination of basic open sets in mind.

For the purposes of the Decoding Lemma \ref{lm:uniqueness} below, the important point is in the other direction, and it states that given a set $A$ for which we know that $A$ is the code for some $H$, there is
only one such $H$ and from $c(H)$ we can read a $\sigma$-Boolean
combination of basic open sets which gives us $H$. This decoding is done in an MSO-definable way, as will be proven in Main Lemma \ref{path:Pi_alpha}.  
\end{example}

\begin{lemma}[{\bf The Decoding Lemma}]\label{lm:uniqueness} Suppose that for some $\alpha\in [1,\omega)$ we are given a set $S$ of branches of $T^{1+\alpha}$. 
Then:
{\noindent (1)} There is an algorithm to determine if
$S=c(H)$ for some set $H$, necessarily a $\Sigma_{\alpha}\cup \Pi_\alpha$-set, and

{\noindent (2)} If there is such $H$, then it is unique. 
\end{lemma}

\begin{proof} The proof if by induction on $\alpha$, proving (1) and (2) simultaneously. 

For \underline{$\alpha=0$}, in the case of closed or open sets, we have that 
$c(H)=H$ with the identification between ${}^{\omega}2$ and $2^\omega$ that is used throughout. We show how to determine if a set $H\subseteq [T]$ is closed. Let $s(H)=\{s\in T:\,(\exists \rho \in H) s\subseteq \rho\}$. This set is a tree and $H$ is closed iff $H=[s(H)]$, in which case $s(H)=T_H$. To determine is $H$ is open, we determine if $2^\omega\setminus H$ is closed. 

For \underline{$\alpha+1$} suppose that $S$ is a set of branches $T^{1+\alpha+1}$. We wish to determine if it is a code for some Borel set $H$ and if
so find find such $H$. It follows from Definition \ref{def:Borelcodes} that any such $H$ must be a $\Sigma_{\alpha+1}$ or a  $\Pi_{\alpha+1}$-set. 
If $S=\emptyset$ then it is the same set as $\emptyset \subseteq T$ and hence we have already determined that it codes the unique Borel set $H=\emptyset$. So let us suppose that $S\neq \emptyset$.

We first determine if $S$ codes an $\Sigma_{\alpha+1}$-set or a $\Pi_{\alpha+1}$-set. 

First we  consider
$\{t(\alpha):\,t\in S\}$. If this is just the constantly 0 branch of $T$, then we know that $S$ might code a $\Sigma_{\alpha+1}$-set, but not any other kind of a Borel set. 
To determine the set possibly coded by $S$, for every $n<\omega$ we consider the following set 
\[
S_n=\{ t\rest [1+\alpha) :\,\pi\in S, t\in \pi, t(\alpha)=(0)_{n}\}.
\]
By the induction hypothesis (1) applied at $\alpha$ we can determine for each $S_n$ if it codes a $\Pi_\alpha$-set $H_n$, and if it does, find a
unique such $H_n$. If any of $S_n$'s fails to code a $\Pi_\alpha$-set, then $S$ does not code a Borel set. Otherwise, suppose we find that $S_n=c(H_n)$ for some
$\Pi_\alpha$-set $H_n$, necessarily unique by the induction hypothesis (2) applied at $\alpha$. If the sequence $\langle H_n:\,n<\omega\rangle$ is not
$\subseteq$-increasing or $H_0=\emptyset$, then $S$ is not a code for a Borel set. Otherwise,
by the definition of codes for $\Sigma_{\alpha+1}$-sets we have that $S=c(H)$ for the set $H=\bigcup_{n<\omega} H_n$.

Now suppose that $\{t(\alpha):\,t\in S\}$ is not a constantly 0 branch of $T$. If it is anything else but the full $T$, we have that $S$ cannot be a Borel code. So suppose that $\{t(\alpha):\,t\in S\}=T$ and let us verify if $S$ codes a $\Pi_{\alpha+1}$-set. For any $n<\omega$ let
$S_n=\{\pi^n:\,\pi\in S\}$. By the induction hypothesis (1) applied at $\alpha$ we can determine for each $S_n$ if it codes a $\Sigma_\alpha$-set $H_n$, and if it does, find such a unique $H_n$. If any $S_n$ fails to code a $\Sigma_\alpha$-set, then $S$ is not a code for a Borel set, so suppose otherwise. 
If the sequence $\langle H_n:\,n<\omega\rangle$ is not
$\subseteq$-decreasing, then $S$ is not a code for a Borel set. Otherwise,
it follows
by the definition of codes for $\Pi_{\alpha+1}$-sets that $S=c(H)$ for $H=\bigcap_{n<\omega} H_n$.

For \underline{$\alpha> 0$} limit, suppose that $S$ is a set of branches of $T^\alpha$. 
Let 
\[
B=\{\beta<\alpha:\,(\exists t\in \bigcup S):\,t(\beta)\neq \Lambda \wedge t\rest (\beta, \alpha)\equiv\Lambda\}.
\]
If $\sup(B)<\alpha$ or $\sup(B)=\alpha$  and the order type of $B$ is not $\omega$, then $S$ does not code a Borel set. Otherwise, we have 
$\sup(B)=\alpha$ and $B$ may be enumerated increasingly as $\langle \beta_n:\,n<\omega\rangle$. (When rereading this proof to make it a proof of Main Lemma \ref{path:Pi_alpha}, we also ask if this is the first such enumeration in a fixed MSO-definable order of countable sequences of countable ordinals and reject $c(H)$ if this is not the case).

Now we ask if it is the case that the set
\[
\{t(\beta_n):\, n<\omega, t\in \bigcup S, t\rest (\beta_n\setminus \alpha)\equiv\Lambda\}\setminus \{\Lambda\}
\]
is a singleton. If it is not, then $S$ is not a Borel code. So suppose that it is a singleton, say $\{l\}$ for some $l<2$. Now for each
$n$ consider the set 
\[
S_n=\{\pi\rest(1+\beta_n):\,\pi\in S \wedge (\forall t\in \pi)[ t(\beta_n)\in \{\Lambda, l\} \wedge t\rest (\beta_n, \alpha)\equiv \Lambda]\}.
\]
By the induction hypothesis, we can determine if $S_n$ is a Borel code and if so find a unique $H_n$ such that $S=c(H_n)$. If any of the $S_n$ fails to be
a Borel code, then $S$ is not a Borel code. So suppose that for every $n$ we have $S_n=c(H_n)$ for some Borel set $H_n$. If $l=0$ we ask if $S$ is the set
of all branches $\pi$ of $T^\alpha$ such that there exists $n$ with $\pi\rest(1+\beta_n)\in S_n$ and $(\forall t\in \pi)[ t(\beta_n)\in \{\Lambda, 0\} \wedge t\rest (\beta_n, \alpha)\equiv \Lambda]$ and if the sets $\langle H_n:\,n<\omega\rangle$ are $\subseteq$-increasing with $H_0\neq \emptyset$. If not, then $S$ is not a Borel code, otherwise
$S=c(H)$ for $H=\bigcup_{n<\omega} H_n$.

Finally, if $l=1$ we ask if $S$ is the set
of all branches $\pi$ of $T^\alpha$ such that there exists $n$ with $\pi\rest(1+\beta_n)\in S_n$ and $(\forall t\in \pi)[ t(\beta_n)\in \{\Lambda, 1\} \wedge t\rest (\beta_n, \alpha)\equiv \Lambda]$, and if the sets $\langle H_n:\,n<\omega\rangle$ are $\subseteq$-decreasing. If so, then $S=c(H)$ for $H=\bigcap_{n<\omega} H_n$. If not, then $S$ is not a Borel code.
\end{proof}

\begin{remark}\label{rem:whyforall} It might seem more intuitive to code $\Sigma_{\alpha+1}$-sets using the existential quantifier, so to take the
set of all branches $\rho^\ast$ in $T^{1+\alpha}\cdot B$ such that for some $n$ the set $\{\rho:\,\rho\frown (0)_n\subseteq \rho^\ast\} $ is a branch in $c(H_n)$.
But that would create problems in the Decoding Lemma at the step of $\Sigma_{\alpha+1}$, since taking all $\rho$ such that for some $\rho^\ast\in c(H)$
we have that $\rho\frown (0)_n\subseteq \rho^\ast$ could give us a set of branches possibly larger than $c(H_n)$.
\end{remark}

\begin{definition}\label{def:eq-codes} Suppose that $1<\alpha<\omega_1$. By $\sim_\alpha$ we denote the equivalence relation on the
set of all subsets of $[T^{1+\alpha}]$ defined by
\[
A\sim_\alpha B 
\begin{cases}
\mbox{if there is a $\Sigma_\alpha\cup \Pi_\alpha$-set $H$ such that both $A$ and $B$ code $H$},\\
\mbox{or neither $A$ nor $B$ code a $\Sigma_\alpha\cup \Pi_\alpha$-set}.
\end{cases}
\]
\end{definition}

By the Decoding Lemma \ref{lm:uniqueness}, each $\sim_\alpha$ is a well-defined equivalence relation on $\PP([T^{1+\alpha}])$. Let $\mathcal C_\alpha$
be the set of all elements of $\PP([T^{1+\alpha}])$ which are codes for $\Sigma_\alpha\cup \Pi_\alpha$-sets and let $\mathfrak C_\alpha=\CC_\alpha/\sim_\alpha$. We shall use ${\mathfrak c}(H)^\alpha$ to denote the $\sim_\alpha$ equivalence class of
$c(H)$, or ${\mathfrak c}(H)$  if $\alpha$ is clear for the context. 

A basic property of codes for Borel sets is given by the following Lemma \ref{le:intersections}. It in particular shows that $\mathfrak C_\alpha$
is
closed under intersection and union. The structure $\mathfrak C_\alpha$ also has a natural order $\subseteq$ and the minimal and maximal element, namely
$\emptyset={\mathfrak c}(\emptyset)^\alpha$ and $\cfr(2^\omega)^\alpha$. It makes it then tempting to say that we are just about to obtain a Boolean algebra, but in fact no, since we do not have a natural notion of complement. That is a consequence of coding by introducing a new copy of $T$ for 
$\Sigma_{\alpha+1}$-sets. This point is crucial in the approach taken here, even though it makes false the intuitive expectation of having a Boolean algebra of codes.

\begin{lemma}\label{le:intersections} (1) Suppose that $\alpha<\omega_1$. Then  if $H_0, H_1$ are $\Sigma_\alpha$-sets ($\Pi_\alpha$-sets respectively),  there are codes $c(H_0)$ and $c(H_1)$ for $H_0$ and $H_1$ respectively such that 
$c(H_0)\cap c(H_1)$ is a code for $H_0\cap H_1$.

{\noindent (2)} The same statement as in (1) is true for $c(H_0)\cup c(H_1)$ representing the codes for $H_0\cup H_1$.
\end{lemma}

\begin{proof} (1) The proof is by induction on $\alpha$. 

In the case \underline{$\alpha=0$}, the conclusion is clear since both open and closed subsets of $2^\omega$ are closed under intersections.

In the case \underline{$\alpha+1$}, suppose that we have two non-empty $\Sigma_{\alpha+1}$-sets $H_l=\bigcup_{n<\omega} H^l_n$ for $l<2$ such that
each $H_l$ is coded using the sequence $\langle c(H^l_n):\,n<\omega\rangle$. In particular, it follows that each sequence 
 $\langle H^l_n:\,n<\omega\rangle$ is $\subseteq$-increasing with $n$ and that $H^l_0\neq \emptyset$.
 
Using De Morgan laws, we have
$H_0\cap H_1= \bigcup_{n<\omega}  \bigcup_{m <\omega} (H^0_n \cap H^1_m)$. We claim that in fact $H_0\cap H_1= \bigcup_{n<\omega} (H^0_n \cap H^1_n)$. To see that, first notice that clearly $\bigcup_{n<\omega} (H^0_n \cap H^1_n) \subseteq  \bigcup_{n<\omega}  \bigcup_{m <\omega} (H^0_n \cap H^1_m)$. In the other direction, we shall use the fact that the representations were chosen as increasing. Namely, suppose that $x\in \bigcup_{n<\omega}  \bigcup_{m <\omega} (H^0_n \cap H^1_m)$. Hence, there is $n$ such that $x\in \bigcup_{m <\omega} (H^0_n \cap H^1_m)$. Since $H^0_n$s are
$\subseteq$-increasing, there is $n_0$ such that for all $n\ge n_0$ we have $x\in \bigcup_{m <\omega} (H^0_n \cap H^1_m)$. Also, $x\in 
\bigcup_{m <\omega} (H^0_{n_0} \cap H^1_m)$ and hence there is $m_0$ such that for all $m\ge m_0$ we have $x\in H^0_{n_0} \cap H^1_m$.
Let $k=\max\{n_0, m_0\}$. Then $x\in H^0_k \cap H^1_k$, so in particular $x\in \bigcup_{n<\omega} (H^0_n \cap H^1_n)$.

We notice that the sequence $\langle H^0_n \cap H^1_n:\,n<\omega\rangle$ is $\subseteq$-increasing. If $H_0\cap H_1=\emptyset$, then 
$H^0_n \cap H^1_n=\emptyset$ for every $n$ and hence a code for the empty set is the empty set and we are done. Otherwise,
it may be that some finitely many $H^0_n \cap H^1_n$ are empty even though $H_0\cap H_1\neq \emptyset$. If so, we throw away those first finitely many
$n$s, hence changing the codes $c(H_0)$ and $c(H_1)$- which is not a problem since we have only said that there are some codes that work.
Now a code $c(H)$ for $H=H_0\cap H_1$ is obtained as the
set of all branches $\rho^\ast$ of $T^{1+\alpha+1}$ where for every $n$ we have that $\{t \in \rho^\ast: \,{\rm ht}(t (\alpha))=n\}=\rho \frown (0)_{n}$ for some element $\rho$ of $c(H^0_n
\cap H^1_n)$. By the induction hypothesis for $\Pi_\alpha$-sets, we have that $c(H^0_n) \cap c(H^1_n)$ is a code for a $H^0_n\cap H^1_n$, which we shall write
as $c(H^0_n\cap H^1_n)=c(H^0_n) \cap c(H^1_n)$. The notation is somewhat abusive since $H^0_n\cap H^1_n$ has many other codes as well, but choosing this one, we have that any branch of $c(H^0_n\cap H^1_n)$ is a branch of $c(H^0_n) \cap c(H^1_n)$. Coming back to which branches $\rho^\ast$ of $T^{1+\alpha+1}$ are in $c(H)$, we obtain that they are those where for every $n$ we have that $\{t \in \rho^\ast: \,{\rm ht}(t (\alpha))=n\}=\rho \frown (0)_{n}$ for some element $\rho$ of
$c(H^0_n) \cap c(H^1_n)$. Hence we obtain $c(H_0\cap H_1)=c(H_0)\cap c(H_1)$.

Now suppose that we have two $\Pi_{\alpha+1}$-sets $H_l=\bigcap_{n<\omega} H^l_n$ for $l<2$ such that
each $H_l$ is coded using the sequence $\langle H^l_n:\,n<\omega\rangle$. In particular, it follows that each sequence 
 $\langle H^l_n:\,n<\omega\rangle$ is $\subseteq$-decreasing with $n$. We have that $H_0\cap H_1= \bigcap_{n,m<\omega}(H^0_m\cap H^1_n)$
 and by the assumptions that the sequences $H^l_n$ are $\subseteq$-decreasing, it follows that $H_0\cap H_1= \bigcap_{n<\omega}(H^0_n\cap H^1_n)$.
 %take k=\min{n,m} and see that $H^0\cap H^1= \bigcap_{k<\omega}(H^0_k\cap H^1_k)$
 By the inductive assumption applied to $\Sigma_\alpha$-sets we can write $c(H^0_n\cap H^1_n)=c(H^0_n) \cap c(H^1_n)$, with the same abuse of notation as in the case of $\Pi_\alpha$-sets. Then we conclude as in the case of $\Sigma_{\alpha+1}$-sets that $c(H_0\cap H_1)=c(H_0)\cap c(H_1)$.
 \begin{comment}
 We check with branches of $T^{1+\alpha+1}$ were chosen to be in $c(H)$. For each $n,l$ we find $S^l_n\subseteq T^{1+\alpha}$ such that $c(H^l_n)$
 is a set of branches of $S^l_n$. Let $S_n=S^0_n\cap S^1_n$, hence $c(H)=c(H^0_n)\cap c(H^1_n)$ is a set of branches of $S_n$. Then a branch $\pi$
 of $T^{1+\alpha+1}$ is in $c(H)$ if for every $n$ we have that $\pi^n$ is a branch of $c(H)=c(H^0_n)\cap c(H^1_n)$, hence iff $\pi$ is a branch of
 $c(H^0)\cap c(H^1)$.
\end{comment}

For \underline{$\alpha>0$ a limit}, suppose that we have two $\Sigma_{\alpha}$-sets $H_l=\bigcup_{n<\omega} H^l_n$ for $l<2$ such that
each $H^l$ is coded using the sequence $\langle H^l_n:\,n<\omega\rangle$. Hence, each of the sets $H_n^l$ is a $\Pi_{\beta^l_n}$-set for some 
$\beta^l_n<\alpha$, the sequences $\langle H^l_n:\,n<\omega\rangle$ are $\subseteq$-increasing and $\sup_{n<\omega}\beta^l_n=\alpha$. Moreover, the sequence $\langle \beta^l_n:\,n<\omega\rangle$ does not depend on $l$ and is the canonically chosen sequence $\langle \beta_n:\,n<\omega\rangle$ used to define codes in $[T^\alpha]$.

If both of the sets $H_0, H_1$ are $\bigcup_{\beta<\alpha}(\Sigma_\beta\cup \Pi_\beta)$-sets the situation is already covered by the induction hypothesis. Otherwise, at least one of them is a $\Sigma_\alpha$-set without being a
$\bigcup_{\beta<\alpha}(\Sigma_\beta\cup \Pi_\beta)$-set. We can assume by symmetry that it is $H_0$. If $H_1=\emptyset$ then the situation is trivial, so let us assume that $H_1\neq \emptyset$. First suppose that $H_0\cap H_1\neq \emptyset$. Let $a\in H_0\cap H_1$ be any element (for definability purposes, we can choose it to be the first in the lexicographic order of
$[T]$). 

Let us suppose first that $H_1$ is a $\bigcup_{\beta<\alpha}(\Sigma_\beta\cup \Pi_\beta)$-set, so, remembering that each $\Sigma_\beta$-set is a
$\Pi_{\beta+1}$-set, we can assume that $H_1$ is a $\Pi_{\beta_m}$-set for some $m<\omega$. We have that $H_0\cap H_1=\bigcup_{n<\omega} (H^0_n \cap H_1)$. We change the
representation used for $H_1$ to let $H^1_n=H_1$ for all $n\ge m$ and $H^1_n=\{a\}$ for $n<m$.

By the induction hypothesis applied to $\beta_n$, we can write $c(H^0_n \cap H_1)=c(H^0_n)\cap c(H_1)$, for every $n$. 
Then we consider which branches $\pi$ of $T^\alpha$ we use in constructing $c(H)$ in the definition of the Borel codes, and we verify that we obtain
exactly $c(H_0 \cap H_1)=c(H_0)\cap c(H_1)$.

Now suppose that neither $H_0$ nor $H_1$ are $\bigcup_{\beta<\alpha}(\Sigma_\beta\cup \Pi_\beta)$-sets.
We have that $H_0\cap H_1=\bigcup_{n<\omega} (H^0_n\cap H^1_n)$, as in the case of $\Sigma_{\alpha+1}$-sets above.
%De Morgan laws as above.
Then by the induction hypothesis applied to $\beta_n$, for each $n$ we can write, with the same conventions as before $c(H^0_n\cap H^1_n)=c(H^0_n)\cap c(H^1_n)$. Similarly as in the above, we can handle the possibility of the first finitely many $n$ such that $H^0_n\cap H^1_n=\emptyset$, so let us assume that this has been done. At the end we consider which branches $\pi$ of $T^\alpha$ we use in constructing $c(H)$ in the definition of the Borel codes, and we verify that we obtain
exactly $c(H_0 \cap H_1)=c(H_0)\cap c(H_1)$.

It remains to deal with the case $H_0\cap H_1=\emptyset$. This implies that there is $n_0$ such that for all $n\ge n_0$ we have that
$H^0_n\cap H^1_n=\emptyset$ and hence by the induction hypothesis we can assume that $c(H^0_{n_0})\cap c(H^1_{n_0})=\emptyset$. By the construction
of codes at limit ordinals it follows that $c(H_0)\cap c(H_1)=\emptyset$.

The case of $\Pi_\alpha$-sets is exactly the same.
%the only difference in the coding is the dummy last non-trivial coordinate

{\noindent (2)} We indicate the changes that are needed in the proof of (1) to obtain the proof for (2). The case $\alpha=0$ is the same. For the case 
$\alpha +1$ now the $\Sigma_{\alpha+1}$-sets no longer require the use of De Morgan laws, but the $\Pi_{\alpha}+1$-sets do, with respect to $\cup$.
The case of $\alpha>0$ limit remains the same, where we replace $\cap$ by $\cup$.
\end{proof}

\section{Decidability of $T^{1+\alpha}$}\label{sec:automates} 
Towards the goal of answering in the affirmative the question of the decidability of the monadic second order logic on the reals where the monadic quantifier ranges over the Borel subsets of the reals, we develop various translations into sentences of $S2S$. The latter we know decidable by \cite{Rabin}. The main translation is introduced in Definition
\ref{def:S2S}. 

\subsection{The recursive correspondance between MSO formulas over $T$ and Rabin automata}\label{subs:Rabin-main}
We recall the main theorem of \cite{Rabin} and the definitions needed to understand its statement.

\begin{definition}\label{def:S2S}
By {\em an MSO formula} we shall mean an MSO-formula over the structure ${\mathfrak N}_2$, which
stands for the tree $T$ with two predicates $r_0$ and $r_1$ for the concatenation of $0$ and $1$, along with the orders $\le$ and $\le_{\rm lex}$.
\footnote{We follow the well chosen notation of \cite{Rabin},
with small exceptions made available by the advantages of LateX over the typing capabilities in 1969. So, Rabin's $\preceq$ is $\le_{\rm lex}$ here. We also denote the set of labels of a labelled tree by $L$ since the letters $\Sigma$ and $\Sigma_n$ used in  \cite{Rabin} are too close to the notation for Borel sets that is important for us.} By {\em
MSO-definable} we shall mean MSO-definable over this structure. The complete MSO theory of ${\mathfrak N}_2$ is denoted by $S2S$.
\end{definition}

It was shown in \cite[pg.8]{Rabin} that both $\le$ and $\le_{\rm lex}$ are MSO-definable using just
$r_0$ and $r_1$. It follows that the MSO theory of the ordered Cantor tree along with the lexicographic order is a sub-theory of $S2S$. 
The main theorem of \cite{Rabin} is Theorem 1.7, which shows that to every formula of $S2S$ there is an associated Rabin automaton which represents the formula. To recall, this means that there is an elementary recursive procedure which to every such formula $F({\bf A}_1, \ldots {\bf A}_n)$, where ${\bf A}_1, \ldots {\bf A}_n$ are set variables, associates a Rabin automaton ${\mathfrak A}_F$ which accepts exactly the set of all {\em representations} of tuples $(A_1, \ldots A _n)$ such that ${\mathfrak N}_2\satisfies F(A_1, \ldots A _n)$.
Let us recall what is meant by the Rabin automaton and representations (see \cite[pg. 6]{Rabin}).

An automaton ${\mathfrak A}$ is an object consisting of finitely many states forming the set $S$ (its ``hardware'') of which some are designated as initial states, forming
the set $S_0$, and
some are designated as final stages, forming the set $F^*$. Rabin automata have an associated set of labels, denoted $L$, which for the example of $n$-place formulas as above will be the set $2^n$. The ``software'' of the automaton or its ``programme'' is 
a table of moves, which is a function $M$ which to every pair consisting of a state and a label, associates a subset of $S\times S$. {\em Runs} of
a Rabin automaton are defined over trees $T_s=\{t \in T:\,s\le t\}$ (so the basic open sets in $T$) labelled by $L$, so pairs of the form $(T_s, v)$ where $v$ is
some function $v:T_s\to L$.

What is a run is defined inductively. For a labelled tree $(T_s, v)$ 
a run $r$ associates to each $t\in T_s$ a state $r(t)$ in a way that for every $t\in T_s$, the pair $(r(t\frown 0) ,r(t\frown 1))$ is a pair of states in $M (r(t), v(t))$. 
We say that the automaton  ${\mathfrak A}$ {\rm accepts} a labelled tree $(T_s, v)$ if there is some run of ${\mathfrak A}$ such that $r(s)\in S_0$ and
for every branch $\rho$ of $T_s$, the run $r$ on $\rho$ has visited every state in $F^*$ infinitely often [the latter is what is known as the {\em acceptance condition of the Rabin automata}, which differentiates these automata from others on the same set of states and with the same table of moves, but with a different acceptance condition]. In other words, for every $a\in F^*$,
the set $\{n<\omega:\,r(\rho\rest n)=a\}$ is infinite.

It remains to define what we mean by representations of tuples $(A_1, \ldots A _n)$ of subsets of $T$. This is a correspondence which to every such tuple of subsets
associates a labelled tree. Namely, with any $n$-tuple $(A_1, \ldots A _n)$ of subsets of $T$ we associate a copy of $T$ labelled by $2^n$, such that every $s\in T$ is labelled by the $n$-tuple $(\chi_{A_1}(s), \ldots \chi_{A_n}(s))$. This is clearly a bijection $\tau$ between $\PP(T)^n$
and the set of labelled trees of the form $(T, v)$ where $v$ is a labelling with values in $2^n$. Therefore, keeping $\tau$ in mind, we can talk about runs of an automaton as accepting or not accepting an $n$-tuple $(A_1, \ldots A _n)$ of subsets of $T$.

Then the main Theorem 1.7 of \cite{Rabin} gives an elementary recursive procedure which to a given a formula $F({\bf A}_1, \ldots {\bf A}_n)$  associates a 
Rabin automaton ${\mathfrak A}_F$ which accepts exactly the $n$-tuples $(A_1, \ldots A _n)$ of subsets of $T$ such that ${\mathfrak N}_2\satisfies F(A_1, \ldots A _n)$.

\subsection{Towards the automata associated to $n$-tuples of subsets of $T^{1+\alpha}$}
We shall be interested in an analogue of Rabin's theorem for subsets of $T^{1+\alpha}$ for $\alpha<\omega_1$, since we have seen that we have
a coding procedure that represents $\Sigma_\alpha\cup \Pi_\alpha$-subsets of $[T]$ as sets of branches in $T^{1+\alpha}$, while a main technique for
Rabin's decidability proof is to use the fact that the closed subsets of $[T]$ can be represented as sets of branches in $T$.
A natural first try would be to define automata that generalise the Rabin automata described in \S\ref{subs:Rabin-main} by replacing $T$
by $T^{1+\alpha}$ in his definitions. However, we use a different approach. In \cite[\S1.9]{Rabin}, Rabin states: ``The binary tree $T_2=T$, in a certain sense, contains as subtrees all trees with countable branching. For this reason, the decidability of
S2S implies decidability of (monadic) second-order theory of more complicated particular
trees and classes of trees with countable branching. Here we shall treat only the
case $T_\omega$.'' He goes on to prove \cite[Th 1.12]{Rabin}, which shows that S$\omega$S is decidable, where this denotes the MSO theory
${}^{<\omega}\omega$ with the naturally defined successor relations. The proof consists of finding an isomorphic copy of ${}^{<\omega}\omega$ within $(T, \le_{\rm lex})$ and MSO-defining the successor relations, therefore interpreting S$\omega$S as a sub-theory of S2S. 

Our main auxiliary theorem is to show that the theory of every $T^{1+\alpha}$ with the set predicates ranging over $\Pi_\alpha$-sets is decidable, and also the theory of $T^\alpha\cdot B$ with the set predicates ranging over $\Sigma_{\alpha+1}$-sets. This has to be made 
precise by specifying which operations are allowed, which we do below, but the idea is to do for these theories what Rabin has done for $S\omega S$.

\subsection{Translating into $S2S$}\label{subsec:translation}
Let $\alpha<\omega_1$. We let $\le_\alpha$ stand for the order on $T¨^{1+\alpha}$ and let $<_\alpha$ stand for the
strict order version of $\le_\alpha$. 

The following lemma summarises what we need in order to support the argument of Main Lemma \ref{path:Pi_alpha}.

\begin{lemma}\label{lem:orders}  The following hold for any $1\le \alpha<\omega_1$:
\begin{itemize}
\item[(1)] there is an MSO-definable function $f^\alpha:\,T^\alpha\to T$ such that 
\[
\{(f^\alpha (x),f^\alpha(y)):\,x\le_\alpha y \}
\]
is MSO-definable,
\item[(2)]  if $x\in \bigcup_{\beta<\alpha}T^\beta$ then, denoting $\lg(x)$ as the unique $\beta$ such that $x\in T^\beta$, we have that the
function $x\mapsto \lg(x)$ is MSO-definable,
\item[(3)]  the projections $x\mapsto b(x)$ from $T^{1+\alpha+1}\to T$ and $x\mapsto p(x)$ from 
$T^{1+\alpha+1}\to T^{1+\alpha}$ (see Notation \ref{note:projection-base}) are MSO-definable,
\item[(4)]  suppose that $\langle F_\beta({\mathbf A}_1, \ldots {\mathbf A}_n):\,\beta<\alpha\rangle$ is a sequence of MSO-formulas.
Then the function $\beta\mapsto F_\beta({\mathbf A}_1, \ldots {\mathbf A}_n)$ is MSO-definable.
\end{itemize}
\end{lemma}

The proof is inspired by Rabin's proof of \cite[Th 1.12]{Rabin} mentioned above.

\begin{proof} (1) The proof is by induction on $\alpha$. The case \underline{$\alpha=1$} is obvious, since we can just let $f^1$ be the identity function and
have $\le_1=\le_T$.

Case \underline{$\alpha+1$}. By the induction hypothesis there is an MSO-formula $F_\alpha(s,t)$ such that for $s,t\in T$ we have
\[
F_\alpha (s,t) \iff (\exists x, y \in T^\alpha) [s= f^\alpha (x) \wedge  t= f^\alpha (y) \wedge x\le_\alpha y].
\]
Let $g=g^{\alpha+1}:\,T^{\alpha+1}\to T\times T$ be given by $(\bar{x},s)\mapsto (f^\alpha(\bar{x}), s)$, hence by the assumption on $f^\alpha$, the function  $g^{\alpha+1}$ is
MSO-definable. 
We shall first define the relation $R_\alpha\subseteq T\times T$ by letting 
\[
R_\alpha=\{((s_0, s_1), (t_0, t_1)):\, [F_\alpha (s_0, t_0) \wedge s_1=t_1 ]\vee s_1 <_T t_1\}.
\]
Let ${\rm pr}: \, T\times T\to T$ be a recursive bijection and for $s,t\in T$ define
\[
F_{\alpha+1}(s,t)\equiv (\exists u_0, u_1, v_0, v_1) [{\rm pr}(u_0, u_1)\!=s\wedge {\rm pr}(v_0, v_1)\!=t \wedge ((u_0, u_1), (v_0, v_1))\,\in R_\alpha].
\]
Unravelling the definitions, we see that it suffices to let $f^{\alpha+1}={\rm pr}\circ g^{\alpha+1}$. Namely,
\[
\begin{split}
\{(f^{\alpha+1}(x), f^{\alpha+1}(y)):\,x\le_{\alpha+1} y\}=\{({\rm pr}( g^{\alpha+1}(x)), ({\rm pr}( g^{\alpha+1})(y)) :\,x\le_{\alpha+1} y\}= \\
\{{\rm pr}(f^\alpha(x\rest\alpha, x(\alpha)), {\rm pr}(f^\alpha(y\rest\alpha, y(\alpha)) :\,x\le_{\alpha+1} y\}=\\
\{{\rm pr}(f^\alpha(x\rest\alpha, x(\alpha)), {\rm pr}(f^\alpha(y\rest\alpha, y(\alpha)) :\,[x\rest\alpha\le_\alpha y\rest\alpha \wedge x(\alpha)=y(\alpha) ]\vee 
x(\alpha)<_Ty(\alpha)\} =\\
\{{\rm pr}(f^\alpha(x\rest\alpha, x(\alpha)), {\rm pr}(f^\alpha(y\rest\alpha, y(\alpha)) :\,F_\alpha(f^\alpha(x), f^\alpha(y)) \wedge x(\alpha)=y(\alpha) ]\vee 
x(\alpha)<_Ty(\alpha)\}=\\
\{{\rm pr}(f^\alpha(x\rest\alpha, x(\alpha)), {\rm pr}(f^\alpha(y\rest\alpha, y(\alpha)) :\,((f^\alpha(x),x(\alpha)),(f^\alpha(y), y(\alpha))\in R_\alpha \}.\\
\end{split}
\]
Therefore for $s,t\in T$ we have 
\[
F_{\alpha+1}(s,t)\iff (\exists x,y\in T^{\alpha+1})( x\le_{\alpha+1} y \wedge s= f^{\alpha+1}(x) \wedge t= f^{\alpha+1}(y)),
\]
as required.

Case \underline{$\alpha>0$ a limit ordinal}. Let us first work with $\alpha=\omega$. The way that Rabin showed in \cite[Th 1.12]{Rabin} that we can MSO-interpret $S\omega S$  within $S2S$ is by finding an MSO-definable
subset $A^*$ of $T$ such that
\begin{itemize}
\item $A^\ast$ is a subtree of $T$ of height $\omega$
\item 
for any $a\in A^\ast$, the set $\{b\in A^*: \,a\le_T b \wedge (a,b)\cap A^\ast=\emptyset \}$ ordered by $<_{\rm lex}$ is order-isomorphic to 
$\omega$.
\end{itemize}
Hence to prove our claim for $\alpha=\omega$ it suffices to show that
$(T^\omega, \le_\omega)$ is
MSO-definable over $S\omega S$. We consider the subset $A$ of ${}^{<\omega}\omega$ of all sequences $\sigma$ satisfying the
following conditions:
\begin{itemize}
\item $\sigma\in {}^{<\omega}\omega$,
\item $\sigma$ can be written as a finite sequence of sequences of elements of $T$, separated by 2, so
\[
\bar{s_0}2\bar{s_1}2 \ldots \bar{s_n}2
\]
for some $n<\omega$. We denote $n=n(\sigma)$.
\end{itemize}
We order $\sigma\le_A \sigma'$ if $n(\sigma)=n(\sigma')$ and for every $i<n$ we have $s_i\le_{T^{\lg(s_i)}} s_i'$.
Then we have that $\le_A$ is MSO-definable and that $(A, \le_A)$ is order-isomorphic to $T^\omega$. Hence, using the induction hypothesis that
every $(T^n,\le_n)$ is MSO-definable, we have given an MSO-definition (over $S\omega S$) of $(T^\omega, \le_\omega)$.

We would like to generalise this to any given countable limit ordinal $\alpha$, but we need a different approach since it is not obvious how to MSO-define
``being of order type $\alpha$''. To see that this possible, we 
first recall B\"uchi's result from \cite{Buchi-ordinals} which says that for any
$\beta<\omega_1$ to any MSO-formula $F(\bar{X})$ over the structure $(\beta,<)$, there corresponds a B\"uchi automaton $\mathfrak A$ which accepts
the input $\bar{X}$ iff $F(\bar{X})$ holds, and this correspondence is also true in the other sense. Then we recall Rabin's result in \cite[\S3]{Rabin} that to the B\"uchi automaton 
$\frak A$ there corresponds a Rabin automaton $\mathfrak B$ and a recursive transformation $f$ of the input space, such that $\mathfrak B$ accepts the input $f(\bar{X})$ iff $\mathfrak A$ accepts
the input $\bar{X}$. (This result uses Mc Naughton's theorem from \cite{zbMATH03336827} which transforms B\"uchi automata to Müller ones). Finally, Rabin's result from \cite{Rabin} gives that that there is an MSO-formula $F_T(\bar{Y})$ such that $F_T(\bar{Y})$ is true over $T$ iff $\mathfrak B$ accepts the input $(\bar{Y})$. Putting these results together, we obtain the interpretation of the MSO theory of $(\beta,<)$ in the
MSO theory of $(T, \le, \le_{\rm lex})$.

Now we have that having the order type $\beta$ is MSO-definable over $(\beta, <)$ by the formula $F(X)=(\forall x)(x\in X)$.
Therefore there is a formula $F_\beta(X, Y)$ such that $F_\beta(X, Y)$ holds in $T$ iff $X\subseteq T$, $Y$ is an order relation on $X$ and
$X$ has order type $\beta$ in the relation $Y$.

Coming back to our $\alpha$, we use these observations to give an MSO-representation of $(T^\alpha,\le_\alpha)$ within $T$. By the
induction hypothesis, we have an MSO-representation of each $(T^\beta,\le_\beta)$ for $\beta<\alpha$ AND the elements of $T^\beta$ are represented by
certain $\beta$-sequences of the elements of $T$. Let $(A_\alpha, \prec_\alpha)$ be an MSO-definable subset of $T$ ordered by $\prec_\alpha$ in the order type
$\alpha$. We shall define a formula $G_\alpha(\varrho)$:
\begin{equation*}
\begin{split}
G_\alpha(\varrho)\equiv \varrho\subseteq A_\alpha\times T \wedge (\forall x, y, z) [(x,y), (x,z)\in \varrho\implies y=z]\wedge \\
(\forall x\in A_\alpha)(\exists y)(x,y)\in \varrho \wedge 
(\exists x\in A_\alpha)(\forall y,z)[(x\prec_\alpha y) \wedge (y,z)\in \varrho \implies z=\Lambda].
\end{split}
\end{equation*}
Hence $G_\alpha(\varrho)$ holds if $\varrho$ is an $\alpha$-sequence of elements of $T$ which is eventually $\Lambda$.
It follows that there is a recursive bijection between $I_\alpha=\{\varrho:\,G_\alpha(\varrho)\}$ and  $T^\alpha$. Hence it suffices to define $f^\alpha$ on $I_\alpha$. For every $\varrho\in I_\alpha$ and $\beta<\alpha$
we can MSO-define $\varrho\rest\beta$ with parameter $\varrho$, since $\beta$ is MSO-definable, ``being an ordinal'' is MSO-definable and ``being an ordinal
$<\alpha$'' is  MSO-definable. We can also MSO-define $\beta(\varrho)$ as the first such that for all $\beta>\beta(\varrho)$ we have $\varrho(\beta)=\lambda$.
Then it suffices to let $f^\alpha(\varrho)=f^{\beta(\varrho)}(\varrho)$. We then have 
\[
\varrho\le_\alpha '\varrho'\iff (\beta(\varrho)= \beta(\varrho') \wedge \varrho\rest\beta(\varrho) \le_{\beta{\varrho}}\varrho'\rest\beta(\varrho)).
\]
This is an MSO-definition of $\le_\alpha$, as claimed to exist.

\smallskip 
{(2)} As in (1), we use the fact that being of order type $(\beta,<)$ is MSO-definable . Given $\bar{x}\in \bigcup_{\beta<\alpha} T^\beta$ we treat $\bar{x}$ as
a function and define the order on ${\rm rg}(\bar{x})$ so that $i<j\implies x(i)<x(j)$. Then for any $\beta<\alpha$ we can decide if ${\rm rg}(\bar{x})$ is of order type 
$\beta$. We let $\beta(\gamma)$ be the unique $\beta$ such that the order type of ${\rm rg}(\bar{x})$ is $\beta$.

\smallskip 
{(3)} Given that we have established that being in $T^{1+\alpha}$ is MSO-definable, this item is similar to the classical case: $b(\bar{x})=y$ if there is $\bar{z}$ such that 
$\bar{x}=\bar{z}\frown y$ and $\bar{z}$ is in $T^{1+\alpha}$. Similarly for the projection.

\smallskip 
{(4)} Similar to the last part of the proof of (2).
\end{proof}

\subsection{Translation between MSO-formulas over $T^{1+\alpha}$ and those over $T$}\label{subsec:no-matter-alpha}
Thanks to Lemma \ref{lem:orders}, we can speak of {\em MSO-definable subsets of} $T^{1+\alpha}$, where a subset $A$ of $T^{1+\alpha}$ is MSO-definable iff the image
${f^{1+\alpha}}(A)=\{f^{1+\alpha}(a):\,a\in A\}$ is MSO-definable as a subset of $T$. 

\begin{notation}\label{not:Malpha} We write ${\mathfrak M}_\alpha$ for the structure $(T^{1+\alpha}, \le_{1+\alpha})$, for $1<\alpha<\omega_1$.
\end{notation}

We can define a {\em translation}
of MSO-formulas over ${\mathfrak M}_\alpha$ to MSO-formulas over ${\mathfrak N}_2$, as follows.

Let $F$ be an MSO-formula $F(\Aa_1, \ldots, \Aa_n)$ with the variables ranging over  $T^{1+\alpha}$. In this notation we use the usual convention that the individual variable $x$ is identified with the set variable $\{x\}$. Then $(tF)=(tF)_{1+\alpha}$ is the MSO-formula
\[
\begin{split}
(tF)({\mathbf C_1},  \ldots, {\mathbf C}_n)\iff (\exists \Aa_1)\ldots (\exists \Aa_n) [(\forall x) [x\in \bigcup_{i=1}^n\Aa_i \implies \lg(x)=1+\alpha] \wedge \\
\bigwedge_{i=1}^n [ {\mathbf C}_i=f^{1+\alpha}(\Aa_i)]\wedge F(\Aa_1, \ldots, \Aa_n)].
\end{split}
\]
Note that we have used the fact that ${\mathbf C}=f^{1+\alpha}(\Aa)$ is an MSO-formula given by
\[
\mathbf C=f^{1+\alpha}(\Aa)\equiv (\forall y)[y\in \mathbf C\iff (\exists x)[x\in \Aa\wedge y= f^{1+\alpha}(x)].
\]
We conclude that for every $F(\Aa_1, \ldots, \Aa_n)$,
\[
{\mathfrak M}_{1+\alpha}\satisfies F(\Aa_1, \ldots, \Aa_n) \iff {\mathfrak N}_2\satisfies (tF)_{1+\alpha}(f^{1+\alpha}(\Aa_1), \ldots, f^{1+\alpha}(\Aa_n)) .
\]

\begin{corollary}\label{cor:dec-alpha} For any $\alpha<\omega_1$, the MSO theory of ${\mathfrak M}_{1+\alpha}$-is decidable.
\end{corollary}

\subsection{Subsets of $T^{1+\alpha}$ and branches in $[T^{1+\alpha}]$} In \cite[\S 2.3]{Rabin}, Rabin showed how to consider MSO-formulas over the real order where only closed and $F_\sigma$-subsets of $2^\omega$ are allowed for the set variables, as being MSO-formulas over $T$. This translation, along with the decidability 
of S2S that he proved in his Main Theorem \cite[Th 1.1]{Rabin}  allowed him to show in  \cite[Th 2.17]{Rabin} that the MSO theory of the real order where only closed and $F_\sigma$-subsets of $2^\omega$ are allowed for the set variables, is decidable. Our main Theorem \ref{th:main} states that the MSO theory of the real order where Borel sets are allowed for the set variables, is decidable. To prove it, we follow the reasoning of Rabin's proof but we use subsets of
$T^{1+\alpha}$ to represent $\Sigma_\alpha\cup\Pi_\alpha$-sets through the coding we developed. The most important remaining part of the translation 
is done through the Main Lemma \ref{path:Pi_alpha}. Before embarking to it, we develop a technique which will let us go back and forth between 
subsets of $T^{1+\alpha}$ and sets of branches in $T^{1+\alpha}$ in an MSO-definable way. This is the same technique that Rabin used for $T$ in 
\cite[\S 2.3]{Rabin}, only we apply it to $T^{1+\alpha}$. Let us develop it now. 

Let us first remark that for every $\alpha$ there is a formula ${\rm Path}_{1+\alpha}(\Bb)$ such that for $B\subseteq T^{1+\alpha}$ we have 
${\mathfrak N}_2 \models (t{\rm Path}_{1+\alpha})(f_{1+\alpha} (B))$
iff $B$ is a path of $T^{1+\alpha}$. This formula is
\[
\begin{split}
[(\forall x)(\forall y) (x\in \Bb \wedge y\in \Bb \!\implies\! x\le_\alpha y \vee y\le_\alpha x)] \wedge\\
(\forall \beta<\alpha)\neg(\exists z)(\forall x)[x\in \Bb \!\implies\! x\rest \beta\le z\rest \beta]].
\end{split}
\]
There is also a formula ${\rm Br}_{1+\alpha}(\Bb)$ which denotes that a path is a branch of $T^{1+\alpha}$. It is the conjunction of ${\rm Path}_{1+\alpha}(\Bb)$ with ${\rm Down}_{1+\alpha}(\Bb)$ where the latter is an MSO-definition of being downward closed, namely
\[
{\rm Down}_{1+\alpha}(\Bb)\equiv (\forall x)(\forall y)[(x\in \Bb \wedge y\le_{1+\alpha} x) \implies y \in \Bb].
\]

We now make a correspondance between subsets $A$ of $T^{1+\alpha}$ and sets of branches in $T^{1+\alpha}$. To an MSO-formula $F(\Bb,\Aa)$ over
$T^{1+\alpha}$ and a set $A\subseteq T^{1+\alpha}$ we associate the set $S^{1+\alpha}_F(A)$ of branches in $T¨^{1+\alpha}$ given by
\[
S^{1+\alpha}_F(A)=\{\pi:\,{\rm Br}_{1+\alpha}(\pi) \mbox{ and } {\mathfrak M}_{1+\alpha}\models F(\pi, A)\}.
\]
From now on we shall use the letter $\pi$ to denote a branch, specifying in which $\le_\alpha$ unless clear from the context. The following observation will be 
useful in the proof of Theorem \ref{th:main}.

\begin{observation}\label{obs:subset} Suppose that $F$ is given as above. Then for $A,  A'\subseteq T^{1+\alpha}$ we have
\[
A \subseteq A' \iff S^{1+\alpha}_F(A)\subseteq S^{1+\alpha}_F(A').
\]
\end{observation}

What remains to be done at this moment is to reformulate the Decoding Lemma \ref{lm:uniqueness} in the appropriate terms. 

\begin{lemma}[{\bf Main Lemma}]\label{path:Pi_alpha} Suppose that $\alpha<\omega_1$. Then 

{\noindent (1)} there is an MSO-formulas $F_{\Sigma_\alpha}(\mathbf B, \mathbf A)$over $T^{1+\alpha}$ such that for any subset $A$ of $T^{1+\alpha}$ we have
\[
S^{1+\alpha}_{F_{\Sigma_\alpha}} (A)= c(H)\mbox{ for some }\Sigma_\alpha\mbox{-set }H\subseteq [T],
\]
and vice versa, for any $\Sigma_\alpha$-set $H\subseteq [T]$ and its code $c(H)$, there is $A\subseteq T^{1+\alpha}$ with 
$S^{1+\alpha}_{F_{\Sigma_\alpha}} (A)= c(H)$.

{\noindent (2)} there is an MSO-formulas $F_{\Pi_\alpha}(\mathbf B, \mathbf A)$over $T^{1+\alpha}$ such that for any subset $A$ of $T^{1+\alpha}$ we have
\[
S^{1+\alpha}_{F_{\Pi_\alpha}} (A)= c(H)\mbox{ for some }\Pi_\alpha\mbox{-set }H\subseteq [T],
\]
and vice versa, for any $\Pi_\alpha$-set $H\subseteq [T]$ and its code $c(H)$, there is $A\subseteq T^{1+\alpha}$ with $S^{1+\alpha}_{F_{\Pi_\alpha}} (A)= c(H)$.
\end{lemma}

\begin{proof} Within Definition \ref{def:Borelcodes} we have briefly noted that all the choices made in the defining the codes can be made in an MSO-definable way, and these statements are justified by Lemma \ref{lem:orders}. Consequently, all the decoding done in the Decoding Lemma \ref{lm:uniqueness} is also done in an MSO-definable way. Hence the proof of the present lemma is exactly the proof of the Decoding Lemma \ref{lm:uniqueness}, observing at each stage that we have used  an MSO-definable procedure. In other words, at each step of the algorithm of the Decoding Lemma we have asked a question that can be formulated as the satisfaction of some MSO-formula over the appropriate $T^{1+\beta}$, and hence by Lemma \ref{lem:orders}, is expressible by an MSO-formula over $T$.
\end{proof}

\subsection{Choosing the $\sigma$-representation by open sets in an MSO-definable way}\label{subs:finitely-many}
The main idea of the proof of Main Theorem \ref{th:main} is a translation of an MSO-formula with Borel set quantifiers into an MSO-formula in S2S.
An important point is to handle the formulas of the sort ${\mathbf v}\in {\mathbf A}$ for an individual variable ${\mathbf v}$ and a set variable 
${\mathbf A}$. This is done by applying Lemma \ref{le:intersections} finitely many times, but it has to be done in an MSO-definable way. We prove a Lemma that shows that this is possible.

\begin{notation} The {\em complexity} of a Borel set is the first level of the Borel hierarchy to which the set belongs. For the purpose of the following lemma \ref{le:representations}, we consider basic open sets in $2^\omega$ and their finite Boolean combinations as having the Borel complexity $-1$.
\end{notation}

\begin{lemma}\label{le:representations} Let $\alpha<\omega_1$. Suppose that $H_0$ and $H_1$ are two Borel sets coded using $[T^{1+\alpha}]$, given 
as $\sigma$-combinations  $\langle H^l_n:\,n<\omega\rangle$ of Borel sets of complexity $<\alpha$ and such that $c(H_l)$ can be constructed from these
$\sigma$-combinations, for $l<2$.
Then there is an MSO-definable procedure for changing the representations $\langle H^l_n:\,n<\omega\rangle$ so to obtain ones satisfying that
$c(H_0\cap H_1)=c(H_0)\cap c(H_1)$.
\end{lemma}

\begin{proof} The proof is by induction on $\alpha$.

Suppose \underline{$\alpha=0$}. For a closed set $H$ we consider the tree $T_H=\{s\in T: (\exists \rho\in H) [s]\subseteq \rho\}$
and enumerate that increasingly as $\{s_n:\, n<\omega\}$ in the lexicographic order on $T$. We choose the
representation $H=\bigcap_{n<\omega}\bigcap_{k<n} [2^\omega\setminus [s_k]]$. In other words, we represent $H$ using the fixed sequence
$\langle \bigcap_{k<n} [2^\omega\setminus [s_k]]:\,n<\omega\rangle$ of Borel sets of complexity $-1$ where the $s_n$s are given in the increasing lexicographic order
and then we obtain a decreasing sequence $\langle H_n:\,n<\omega\rangle$ by taking $H_n=\bigcap_{k<n} [2^\omega\setminus [s_k]]$.

For the open sets, we use the complements of what we have just done. Notice that $[\Lambda]=2^\omega=[s_0]$. These representations have the property that $c(H)=H$ for both open and closed sets, so the required intersection property follows.

For the case \underline{$\alpha+1$}, let us first consider the $\Sigma_{\alpha+1}$-sets $H_0$ and $H_1$, represented using the sequences
$\langle H^l_n:\,n<\omega\rangle$ for $l<2$ of $\Pi_\alpha$-sets. Going back to the proof of Lemma \ref{le:intersections}, in the case that $H_0\cap H_1\neq \emptyset$ we have that given any codes  $\langle H^l_n:\,n<\omega\rangle$ for $l<2$, we can in an MSO-definable way, by looking for the first $n$ such that $H^0_n\cap H^1_n\neq \emptyset$, change the codes so that $\langle H^0_n\cap H^1_n:\,n<\omega\rangle$ is a code for $H^0 \cap H^1$. By the inductive hypothesis
for the $\Pi_\alpha$-sets, for every $n<\omega$ we have an MSO-definable way of choosing the codes for $H^0_n$ and $H^1_n$ so that $c(H^0_n\cap H^1_n)=c(H^0_n)\cap c(H^1_n)$, say by some MSO-definable function $F_n$. By Lemma \ref{lem:orders}(4), we can find an MSO-definable function $F$
which lets us perform the whole sequence $\langle F_n:\,n<\omega\rangle$ of the changes simultaneously. Then $F$ gives us the codes for $H^0$ and $H^1$ which are as required. In the case that $H_0\cap H_1=\emptyset$, we have that for each $n$ the intersection $H^0_n\cap H^1_n$ is empty, hence we apply the induction hypothesis for the $\Pi_\alpha$-sets in the same way as in the previous case. 

In the case of $\Pi_{\alpha+1}$-sets we similarly follow the procedure described in the proof of Lemma \ref{le:intersections} and use the induction hypothesis for the $\Sigma_\alpha$-sets.

For the case \underline{$\alpha>0$ limit}, consider again first the $\Sigma_\alpha$-sets. We have chosen a sequence  $\langle \beta:\,n<\omega\rangle$
with $\sup_{n<\omega}\beta_n=\alpha$ in an MSO-definable way.
We have assumed that both $H_0$ and $H_1$ are coded as
subsets of $[T^\alpha]$, so we can pass immediately to the third case considered in the proof of Lemma \ref{le:intersections}. For each $l$ there
is a chosen representation $H_l=\bigcup_{n<\omega} H^l_n$, where each $H^l_n$
is a $\Pi_{\beta_n}$-set. By the induction hypothesis applied to the $\Pi_{\beta_n}$-sets for each $n$, we have an MSO-definable function which lets us change the codes of $H^0_n$ and $H^1_n$ if
necessary to obtain $c(H^0_n\cap H^1_n)=c(H^0_n)\cap c(H^1_n)$. As in the successor case, using Lemma \ref{lem:orders}(4), we can find an MSO-definable function which lets us do these changes simultaneously and hence we obtain codes for $H_0$ and $H_1$ such that $c(H_0\cap H_1)=
c(H_0)\cap c(H_1)$, as required.

The case of $\Pi_\alpha$-sets is the same.
\end{proof}

\section{Decidability of MSO over the real order with Borel Set Quantifiers}\label{sec:Borel}
We have announced that we shall prove the decidability of the MSO theory of the real order with the set variables ranging over the Borel sets. This is the exact wording of Shelah's Conjectures 7A and 7B from \cite{Shelahmonadic}, which are equivalent. We should first explain what the statement means.

Let $\mathfrak C$ be the structure consisting of the Cantor set $[T]$ with its usual lexicographic order $\le$. What we mean by the usual order on $[T]$  is, denoting by $\Delta(\rho, \rho')$ 
the value $\min\{n< \omega:\rho(n)\neq \rho'(n)\}$ such that $\rho(n)\neq \rho'(n)$ for $\rho\neq \rho'\in [T]$, 
 \[  
 \rho\le \rho' \mbox{ iff } \rho(\Delta(\rho, \rho'))<\rho'(\Delta(\rho, \rho')).
 \]
In other words, this is the lexicographic ordering. Then there is a recursive bi-interpretation between the reals with their usual order and $\mathfrak C$, and moreover this bi-interpretation gives rise to a Borel isomorphism. Hence for the decidability of the MSO theory of the real order with the set variables ranging over the Borel sets (Shelah's Conjecture 7A ) it suffices to prove the decidability of the MSO theory of $\mathfrak C$ and set variables ranging over the Borel subsets of 
$[T]$ (Shelah's Conjecture 7B), and vice versa. 

In the sentence above the statement of Conjecture 7A, Shelah motivates the conjectures by mentioning Rabin's \cite[Th 2.7]{Rabin}, which he cites as ``quantification over $F_\sigma$ sets gives a decidable theory''. The precise statement of \cite[Th 2.7]{Rabin} assumes that the language has separate
symbols for set variables which are closed sets and those that are $F_\sigma$-sets. The analogous language for Borel sets can be
understood in two different ways, which turn out to have a different meaning. We shall first explain how to handle the weaker meaning, which is to 
work with the language that has separate symbols for set variables on each level of the Borel hierarchy. This can also be implemented by adding unary relation symbols $\Sigma_\alpha$ and $\Pi_\alpha$ acting on the set variables. We call the resulting interpretation of the Conjectures 7A/7B from
\cite{Shelahmonadic} the {\bf Weak Shelah Conjecture}. This conjecture is also listed as an open question in \cite[pg. 503]{GurevichMSO}. We confirm the Weak Shelah's conjecture in the following Theorem \ref{th:main}. We confirm the stronger version of the conjecture in Corollary \ref{th:strongmain}. That stronger version is explained after the proof of Theorem \ref{th:main}.

\begin{theorem}\label{th:main} Let $\LL$ be the second order language $\{\in, \le, \Sigma_\alpha, \Pi_\alpha (\alpha<\omega_1)\}$ where
$\in$ is a binary relation between the individual variables and the set variables, $\le$ is a binary relation on the individual variables and 
 $\Sigma_\alpha, \Pi_\alpha (\alpha<\omega_1)$ are unary relation symbols on the set variables.
 
Let $\mathfrak C$ be the $\LL$-structure consisting of the Cantor set $[T]$ with $\in$ interpreted as the membership relation, $\le$ interpreted by the lexicographic order on $[T]$, and each $\Pi_\alpha({\mathbf A})$ and the $\Sigma_\alpha({\mathbf A})$
interpreted as ${\mathbf A}$ belonging to the corresponding level of the Borel hierarchy of the Borel subsets of  $[T]$.

Then the MSO theory of $\mathfrak C$ is decidable.
\end{theorem}

\begin{proof} Let $T$ stand for the MSO theory of $\mathfrak C$. Notice that $T$ includes the following sentences:
\begin{itemize}
\item $(\forall \mathbf A) [\Pi_\alpha({\mathbf A})\implies [\Pi_{\alpha+1}({\mathbf A}) \wedge \Sigma_{\alpha+1}({\mathbf A})]]$ and
\item  $(\forall \mathbf A) [\Sigma_\alpha({\mathbf A})\implies [\Pi_{\alpha+1}({\mathbf A}) \wedge \Sigma_{\alpha+1}({\mathbf A})]]$.
\end{itemize}
Let $\varphi$ be a sentence of $T$. 
Let $\alpha<\omega_1$ be large enough so that all set variables
mentioned in $\varphi$ are within $\Sigma_\alpha\cup \Pi_\alpha$-sets. By increasing $\alpha$ if necessary and passing to a sentence equivalent to $\varphi$ modulo $T$, we can assume that all 
of the set variables mentioned in $\varphi$ are $\Sigma_{\alpha}$.

Recall Notation \ref{not:Malpha}. We shall translate $\varphi$ into an MSO-sentence $\varphi^\ast$
of ${\frak M}_{1+\alpha}$. Let $H_1, \ldots H_n$ be the set variables appearing in $\varphi$. We identify individual variables $\rho$ in $\varphi$ with sets 
$\{\rho\}$, which are non-empty closed sets in $[T]$.
In particular they also are $\Sigma_{\alpha}$-sets. Let $\rho_1, \ldots \rho_m$ be all the individual variables appearing in $\varphi$.
By applying Lemma \ref{le:representations} finitely many times if necessary, we
can assume that we have chosen the codes $c(H_i)$ and $c({\rho}_j)$ for $i<n, j<m$ within $[T^{1+\alpha}]$ so that 
\begin{equation}\label{eq:1}
\begin{split}
\\
\rho_j\in H_i  \!  \iff \!  \{\rho_j\}\subseteq H_i  \!  \iff \!   \{\rho_j\}\cap H_i=\{\rho_j\} \! \iff\! c(\{\rho_j\}\cap H_i)=c(\{\rho_j\})\cap c(H_i) \\
\iff c(\{\rho_j\})=c(\{\rho_j\})\cap c(H_i)\iff c(\{\rho_j\})\subseteq c(H_i).
\end{split}
\end{equation}
Recall the notation from Main Lemma \ref{path:Pi_alpha}(1). We shall define a partial order $\preceq$ among those subsets $A$ of
$T^{1+\alpha}$ which satisfy that $S^{1+\alpha}_{F_{\Sigma_{\alpha}}} (A)= c(H)$ for a code $c(H)$ for a singleton set $\{\rho\}\subseteq [T]$. 
Being a set of that form is MSO-definable, as asking if a set is a singleton is an MSO-question and Main Lemma  \ref{path:Pi_alpha}(1)
gives that asking if a set is a code for a $\Sigma_\alpha$-set is an MSO-question. Say that
$A=A(\rho)$ in the case of a positive answer and define $A(\rho)\preceq A(\rho')$ iff $\rho\le \rho'$ in $[T]$.

An application of Main Lemma \ref{path:Pi_alpha}(1) gives us 
subsets $A_1, \ldots  A_n$ and $B_1, \ldots B_m$ of $T^{1+\alpha}$ 
such that 
for each $i\in \{1, \ldots,n\}$ we have that $S^{1+\alpha}_{F_{\Sigma_{\alpha}}} (A_i)= c(H_i)$ and for each  $j\in \{1, \ldots,m\}$ we have that 
$S^{1+\alpha}_{F_{\Sigma_{\alpha}}} (B_j)= c(\{\rho_j\})$. 

We shall use $\varphi$ to define an MSO-sentence $\varphi$ over $T^{1+\alpha}$. This is done by replacing each appearance of
$H_i$ by $A_i$, each appearance of $\rho_j$ by $B_j$, each appearance of $\rho_j\le \rho_k$
by $B_j\preceq B_k$ and by writing $B_j\subseteq A_i$ iff $\rho_j\in H_i$. Notice that the $\subseteq$-relation is definable from the
membership relation present by definition in the MSO-logic.

Recalling Observation \ref{obs:subset} and using (\ref{eq:1}),
we conclude that the way we defined $\varphi^\ast$ is such that 
\[
{\mathfrak C}\satisfies \varphi\iff {\mathfrak M}_{1+\alpha}\satisfies \varphi^\ast.
\]
By Corollary \ref{cor:dec-alpha}, the truth of $\varphi^\ast$ in ${\mathfrak M}_{1+\alpha}$ is decidable.
Hence, the truth of $\varphi$ within ${\mathfrak C}$ is decidable.
\end{proof}

It is possible to formulate a stronger statement than the statement proved in Theorem \ref{th:main} by interpreting Conjecture 7A/7B from
\cite{Shelahmonadic} using a different language.

\smallskip

{\bf Strong Shelah's Conjecture} Let $\LL^\ast$ be the second order language $\{\in, \le, {\mathcal B}\}$ where
$\in$ is a binary relation between the individual variables and the set variables, $\le$ is a binary relation on the individual variables and 
${\mathcal B}$ is a unary relation symbols on the set variables. Let $\mathfrak C^\ast$ be the $\LL$-structure consisting of the Cantor set $[T]$ with $\in$ interpreted as the membership relation, $\le$ interpreted by the lexicographic order on $[T]$, and ${\mathcal B}$ interpreted as being a Borel set.

Then the MSO theory of $\mathfrak C^\ast$ is decidable.

\smallskip

Theorem \ref{th:main} does not answer the Strong Shelah's Conjecture as the set of MSO-sentences possible to express in 
$\mathfrak C^\ast$ is strictly richer than the analogously defined set using $\mathfrak C$. An example of an MSO-sentence 
which is expressible in $\mathfrak C^\ast$ and not in $\mathfrak C$ is the statement of the Borel Determinacy (see \S\ref{sec:Rabin-Borel-determinacy}). 
All known proofs of the Borel Determinacy involve an induction over the $\omega_1$-many levels of the Borel hierarchy and then making the final conclusion using the fact that if $\varphi(\alpha)$ is true for every $\alpha<\omega_1$ then $(\forall \alpha<\omega_1)\varphi(\alpha)$ is true. This is not 
expressible in S2S. However, we can approach the Strong Shelah's Conjecture using the fact that the MSO theory of $(\omega_1, <)$ is decidable, as proven by Büchi in \cite{Buchi1973}, and in a different way by Shelah in \cite[Conclusion 4.10]{Shelahmonadic}. 
A consequence of this result is the following:

\begin{theorem}[{\bf Büchi}]\label{th:Buchi} Suppose that $\langle F_\alpha:\,\alpha<\omega_1\rangle$ is a decidable sequence of decidable sentences, namely, there is a recursive procedure which on an input $\alpha$ outputs the truth value of $F_\alpha$. Then the statement $[(\forall \alpha<\omega_1)F_\alpha]$ is decidable.
\end{theorem}

\begin{proof} We can reformulate the statement as deciding if $\{\alpha<\omega_1:\,\neg F_\alpha\}=\emptyset$. By the assumption, this is an MSO-sentence over $(\omega_1, <)$, and hence decidable by the fact that the MSO theory of $(\omega_1,<)$ is decidable. 
\end{proof}

Putting Theorem \ref{th:main} and Theorem \ref{th:Buchi} together, we obtain:

\begin{corollary}\label{th:strongmain} The Strong Shelah Conjecture is true.
\end{corollary}

\begin{proof} 
Let $\varphi$ be a sentence in the MSO theory of $\mathfrak C^\ast$. The proof of the decidability of 
$\varphi$ is by induction on the
complexity of $\varphi$. 

As the set of decidable sentences is closed under $\neg, \wedge$ and $\vee$, the nontrivial case is the quantifier case. Suppose first that 
$\varphi$ has the prenex normal form $(\forall {\mathbf A})\psi({\mathbf A})$ where $\psi$ does not have any set quantifiers. 
For each $\alpha<\omega_1$ consider the sentence
$F_\alpha\equiv (\forall {\mathbf A})[(\Sigma_\alpha ({\mathbf A})\vee \Pi_\alpha ({\mathbf A})\implies \psi({\mathbf A})]$. 
This is a sentence in the MSO theory of $\mathfrak C$ and so decidable by Theorem \ref{th:main}. Moreover, considering the sequence 
$\langle F_\alpha:\,\alpha<\omega_1\rangle$, it is what we called a decidable sequence of decidable sentences in Theorem \ref{th:Buchi}.

We have that $\varphi$ holds iff $[(\forall \alpha<\omega_1)F_\alpha]$ holds, and by Theorem \ref{th:Buchi}, this is decidable. 

A similar argument applies to MSO-sentences whose prenex forms start with an existential set quantifier or any finite number $n$ of set quantifiers,
followed by a set-quantifier-free formula. We simply need to instead of $\Sigma_\alpha ({\mathbf A})\vee \Pi_\alpha ({\mathbf A})$
in the definition of $F_\alpha$ above use the appropriate finite combinations of $\Sigma_\alpha$s and $\Pi_\beta$s, hence a conjunction of $\le 2n$
statements of the form $\Sigma_\alpha ({\mathbf A})$, $\Pi_\beta({\mathbf A})$ or their negations. We can enumerate all resulting MSO-sentences
over $\mathfrak C$ in the order type $\omega_1$ and then use Theorem \ref{th:Buchi} as in the case of $(\forall {\mathbf A})$.

Hence we conclude that the MSO-theory of $\mathfrak C^\ast$ is decidable, which is exactly the statement of the Strong Shelah Conjecture.
\end{proof}

\section{Effective Borel Determinacy}\label{sec:Rabin-Borel-determinacy}
We recall the definition of the determinacy game for subsets of $2^\omega$. 

\begin{definition}\label{det:detgame} Suppose that $H\subseteq 2^\omega$. We define $G(H)$ as the game in which two players I and II take turns to
construct an element $\rho\in 2^\omega$ by giving the value $\rho(2n)$ for I and $\rho(2n+1)$ for II at the round $n$. I wins iff $\rho\in H$.

A winning strategy $S$ for I is a function defined on $\bigcup_{n<\omega} {}^{2n} 2$ such that for any play $\rho$ of $G(H)$ such that at every step 
$n$ the player I plays $S(\rho\rest 2n)$, the player I wins (this is a recursive definition). We similarly define a winning strategy for II. The set $H$ is determined if one of the players has a winning strategy in $G(H)$.

{\em Borel Determinacy} is the statement that all Borel sets are determined.
\end{definition}

Clearly, a set is determined iff its complement is determined. It is quite easy to see that finite sets are determined. The following is the classical
Gale-Stewart Theorem from \cite{zbMATH03078993}, see \cite[Th 20.1]{Kechris}.

\begin{theorem}[{\bf Gale-Stewart Theorem}]\label{th:Gale-Stewart} Open (and hence closed) sets are determined.
\end{theorem}

Philip Wolfe \cite{zbMATH03113183} proved the determinacy of $F_\sigma$-sets.
In \cite[\S2.5]{Rabin}, Rabin showed that the statement that every $F_\sigma$-set is determined, is expressible as an MSO-sentence over ${\mathfrak N}_2$. He concluded that therefore, Wolfe's theorem is among those decided by the decision procedure for S2S. In particular, knowing that all $F_\sigma$-sets are determined, we can conclude that for every $F_\sigma$-set there is a recursive procedure to decide which player wins the game. This recursive procedure is
given by a run of the Rabin automaton. This is what we may call an {\em effective version of the determinacy for $F_\sigma$-sets}.

Borel Determinacy is a theorem of ZFC, see \cite[Th 20.5]{Kechris}, attributed to Donald Martin. 
Using the coding of arbitrary Borel sets we introduced, Rabin's method of translating MSO-sentences to the runs of the Rabin automata and the fact that we know that Borel Determinacy is true, we obtain an {\em effective Borel's determinacy}. Namely, for every Determinacy Game $G(H)$ where $H$ is a Borel set, there is a run of a Rabin automaton that determines the player that has the winning strategy.  
To summarise the situation:

\begin{corollary}[{\bf Effective Borel Determinacy}] \label{th:effectiveBorel} Determinacy of any given level of the Borel hierarchy is expressible as a sentence in the MSO-theory of $\mathfrak C$. For every Borel set $H$ there is a recursive procedure given as a run of a Rabin automaton which determines if the first or the second player wins in the Determinacy Game $G(H)$. 
\end{corollary}

\begin{proof} On \cite[pg.17]{Rabin} Rabin explains that strategies in the Determinacy Game can be viewed as characteristic functions of
subsets $A$ of $T$. He deals with $F_\sigma$-sets, but for us, let us simply note that he gives
MSO-formulas $W_L({\mathbf A}, {\mathbf D})$ for $L\in \{I,II\}$  over $T$ which express that $\chi_A$ is a winning strategy for Player L in the Determinacy Game for an $F_\sigma$- or a $G_\delta$-set $D$. The formula $W_L({\mathbf A}, {\mathbf D})$ is a conjunction of an MSO-formula $w_L({\mathbf A}, {\mathbf D})$ 
and a separate formula stating that ${\mathbf D}$ is an $F_\sigma$-set or its complement.

It follows from the translations explained in the proof of Theorem \ref{th:main} that for every $\alpha<\omega_1$ there is an MSO-formula over $T$
which expresses that ${\mathbf D}$ is a $\Sigma_\alpha$-set. Say that this formula is $\Sigma_\alpha({\mathbf D})$. Hence the conjunction of 
$\Sigma_\alpha({\mathbf D})$ and $w_L({\mathbf A}, {\mathbf D})$ gives us an MSO-formula  which expresses that $\chi_{\mathbf A}$ is a winning strategy for Player L in the Determinacy Game $G(\mathbf D)$ for a $\Sigma_\alpha$-set $\mathbf D$. Call the two resulting formulas $W_L^\alpha({\mathbf A}, {\mathbf D})$, for $L\in \{I,II\}$. 

Let $H$ be a Borel set and let $\alpha<\omega_1$ be such that $H$ is a $\Sigma_\alpha$-set. 
Then, by Borel Determinacy, exactly one of the sentences $(\exists {\mathbf A})W_I^\alpha({\mathbf A}, H)$ and $(\exists {\mathbf A})W_{II}^\alpha({\mathbf A}, H)$ holds in $\mathfrak C$.
By Theorem \ref{th:main} for each one of these sentences, there is a Rabin automaton that determines if it is true or not. Hence, we have an effective procedure for determining 
the player who has the winning strategy. Moreover, by the definition of $w_L({\mathbf A}, {\mathbf D})$, it follows that in the case of the existence of such a strategy, we have an effective procedure through a Rabin automaton for
determining such a strategy.
\end{proof}

Combining Corollary \ref{th:effectiveBorel} and Theorem \ref{th:Buchi} as we did in the proof of Theorem \ref{th:strongmain}, we obtain the following
Corollary \ref{cor:expressingBD}. Of course, we already know that the statement of the Borel Determinacy is a true theorem, and we have used that fact in the proof of Corollary \ref{th:effectiveBorel}. What Corollary \ref{cor:expressingBD} states is simply that the statement is among those decidable by the decision procedure of the MSO-theory of $\mathfrak C^\ast$.

\begin{corollary}\label{cor:expressingBD} The statement of the Borel Determinacy is a sentence in the MSO-theory of $\mathfrak C^\ast$, and hence
decidable by the decision method of that theory.
\end{corollary}

\section{Closing Remarks}\label{sec:remarks}
We note that our use of the word ``effective'' in the statement of Theorem \ref{th:effectiveBorel} is in the sense that the word is used in mathematical logic, namely that there is a theoretical computing machine which has a run that computes the problem in question. This does not say anything about an actual implementation of an algorithm for this computation. Such an implementation is clearly non-trivial, because Rabin automata work on trees with infinite branches, hence do not have a physical equivalent. Instead of having a physical Rabin automaton, one can construct a simulator, such as the Rabinizer \cite{10.1007/978-3-319-96145-3_30}, which, roughly speaking, to a given programme of the Rabin automaton associates an executable programme in a real programming language. (In fact Rabinizer can start from a formula of the linear temporal logic, translate it into a run of a deterministic Rabin automaton and execute it through a simulation in Java). In this context one can discuss what the computer scientists would understand by ``effectiveness'', which is a polynomial or better bound on the complexity of the simulation programme in Java. We are not aware of any work studying such complexity even for the determinacy for $F_\sigma$-sets.

The elegant proof of the Borel Determinacy given by Kechris in \cite[Th, 20.5]{Kechris} relates to Martin's paper \cite{zbMATH03993574}. Martin mentions in \cite{zbMATH03993574} that it was proposed by Yannis Moschovakis in his book \cite{zbMATH03674071} that such a proof should exist, although it was not found then. Before this proof, some earlier papers on the subject used a different technique. The idea was to associate to each Borel game $G$ a game $G^\ast$ on an open set in a different Polish space, so that if the player I or II has a winning strategy in $G^\ast$, then the player I or II respectively has a winning strategy in $G$. The Borel Determinacy then follows from the Gale-Stewart theorem on the determinacy of Determinacy Gales on open sets in Polish spaces \cite{zbMATH03078993} (this would require proper formulation of such games, more general than Definition \ref{det:detgame}). The idea of translation between $G$ and $G^*$ was introduced by Martin in \cite{zbMATH03342849}, where he claims the proof of Borel Determinacy from the existence of an $\omega_1$-Erd\H{o}s cardinal. His paper \cite{Boreldet} claims Borel determinacy in ZFC, using a better coding. In the two papers the details of the proofs are missing or incorrect. 
With this in mind, we have tried to prove Borel Determinacy in yet another way, to match Martin's original intuition. This would be by playing games on codes for $\Sigma_\alpha \cup \Pi_\alpha$-sets, understanding those as closed games in an appropriate Polish subspace of $[T^{1+\alpha}]$ (see \S\ref{sec:metric}) and then projecting to the ordinary games on closed sets. 
The intended proof is by induction on $\alpha$ and uses a direct proof for $F_\sigma$-determinacy given in \cite[Exercise 20.9]{Kechris} at the non-trivial part of the projection between games in the successor case. This proof would correspond to the original programme of proving Borel Determinacy as intended by Martin in his papers \cite{zbMATH03342849} and \cite{Boreldet}. We have not worked out all of the details of this proof, as from the purely mathematical point of view it seems unnecessary to give a very long combinatorial argument for a theorem whose proof is already known and well written (\cite[Th 20.5]{Kechris}). In fact, we believe that the combinatorial proof is simply a recounting of the decision procedure for the player which has the winning strategy. 
The part of the combinatorial proof that we did write is available in the {\LaTeX} code of this paper on the ArXiv. If an interested person uses it, they should ask for a permission and their work should have a proper mention of the source (this paper). The reason we mention these issues is that we believe that this algorithm might be of interest in the context of implementing it on the simulators of Rabin automata, such as Rabinizer \cite{10.1007/978-3-319-96145-3_30}.

Another remark on our editorial choices is that we present \S\ref{sec:realcoding} separately and before \S\ref{sec:automates} to emphasise the main ideas, for human reading purposes. If writing
a code for a computer program verifying this paper, these two sections would have to be merged into one large inductive definition and proof of the objects within it being MSO-definable done together in a simultaneous induction. 

We finally mention that the proof of Borel Determinacy from \cite[Th. 20.5]{Kechris} or \cite{zbMATH03993574} has indeed been implemented in Lean by Manthe \cite{Borel-det-Lean}. Manthe has also published on ArXiv an outline of an intended proof of Theorem \ref{th:main} \cite{Manthe-sketch} and has presented 
the same in talks, since October 2024. We would be eager to see the details of the proof if it converges, especially since its announcement uses a very different technique than the one presented here.

\bibliographystyle{plainurl}
\bibliography{../bibliomaster}
\end{document}